%% file: root.tex
\documentclass[10pt,a4paper]{article}
\usepackage{sectsty}
\sectionfont{\large}
\subsectionfont{\normalsize}
\subsubsectionfont{\normalsize}
\paragraphfont{\normalsize}
\allsectionsfont{\bfseries\sffamily}
\usepackage{abstract}

\usepackage{caption}
\usepackage{graphicx} 
\usepackage[cmex10]{amsmath} 
\usepackage{amssymb}  
\usepackage[english]{babel}

\usepackage{amsthm}
\usepackage{thmtools}
\usepackage{dsfont} 
\usepackage[usenames,dvipsnames]{xcolor}
\usepackage{hyperref} 
\usepackage{cleveref}

\usepackage{enumitem} \setlist[itemize,1]{leftmargin=\dimexpr 25pt}
\usepackage{cite}
\usepackage{url}
\usepackage{microtype}
\usepackage{BOONDOX-ds}
\usepackage{pgfplots} \pgfplotsset{compat=1.14} \usepgfplotslibrary{colorbrewer}
\usepackage{tikz} 
\usepackage{quoting}
\usepackage{siunitx}
\usepackage{xstring}
\usepackage[T1]{fontenc}
\usepackage{circuitikz}

\makeatletter
\def\@IEEEsectpunct{.\ \,}
\def\paragraph{\@startsection{paragraph}{4}{\z@}{1.5ex plus 1.5ex minus 0.5ex}%
{0ex}{\normalfont\normalsize\itshape}}
\makeatother 

\declaretheorem[style=definition]{theorem}

\declaretheorem[style=definition]{lemma}
\declaretheorem[style=definition,qed=$\vartriangle$]{remark}

\declaretheorem[style=definition,numbered=no,qed=$\lrcorner$]{standing assumption}
\declaretheorem[style=definition,qed=$\lrcorner$]{assumption}
\declaretheorem[style=definition,qed=$\lrcorner$]{definition}
\declaretheorem[style=definition,qed=$\lrcorner$]{problem}

\newcommand {\nn}{\nonumber}
\newcommand{\beq}{\begin{equation}}
\newcommand{\eeq}{\end{equation}}
\newcommand {\bseq}{\begin{subequations}}
\newcommand {\eseq}{\end{subequations}}
\newcommand {\bma}{\left[}
\newcommand {\ema}{\right]}

\newcommand {\Zplus}{\mathbb{Z}_{+}} 	
\newcommand {\R}{\mathbb{R}} 	
\newcommand {\Rplus}{\mathbb{R}_{+}} 	
\newcommand {\Co}{\mathbb{C}} 	
 	
\newcommand {\Cominus}{\mathbb{C}_{-}} 	
\newcommand {\Cozero}{\mathbb{C}_{0}} 	

\newcommand{\Image}{\operatorname{im}} 
\newcommand{\rank}{\operatorname{rank}} 
\newcommand{\transpose}{\mathsf{T}} 
 
\newcommand{\norm}[1]{\left\lVert#1\right\rVert}
\newcommand{\Span}{\operatorname{span}}
\newcommand{\diag}{\operatorname{diag}}
\newcommand{\spectrum}[1]{{\sigma({#1})}}
\newcommand{\spectrumk}[1]{{\sigma_k({#1})}}

\title{\Large {\sffamily\bfseries  Model reduction by least squares moment matching\,for\,linear\,and\,nonlinear\,systems}}

\author{Alberto Padoan}

\date{}

\begin{document}

\maketitle
\thispagestyle{plain}
\pagestyle{plain}

\begin{abstract}     
\noindent 
The paper addresses the model reduction problem for linear and nonlinear systems using the notion of least squares moment matching. For linear systems, the main idea is to approximate a transfer function by ensuring that the interpolation conditions imposed by moment matching are satisfied in a least squares sense.  
The paper revisits this idea using tools from output regulation theory to provide a new time-domain characterization of least squares moment matching. It is shown that least squares moment matching can be characterized in terms of an optimization problem involving an invariance equation and in terms of the steady-state behavior of an error system. This characterization, in turn, is then used to define a nonlinear enhancement of the notion of least squares moment matching and to develop a model reduction theory for nonlinear systems based on the notion of least squares moment matching. Parameterized families of models achieving least squares moment matching are determined both for linear and nonlinear systems. The new parameterizations are shown to admit natural geometric and system-theoretic interpretations.  The theory is illustrated by worked-out numerical examples. 
\end{abstract}

\section{Introduction}

Model reduction is the art of approximating the behavior of a dynamical system while preserving its main features~\cite{antoulas2005approximation}. This task occurs frequently in control engineering practice and can be posed,  mathematically, as an approximation problem. 
For linear systems, a popular way to approach the model reduction problem is through methods based on moment matching~\cite{antoulas2005approximation,antoulas2020interpolatory,antoulas1990solution,
grimme1997krylov,gallivan2004sylvester,gallivan2006model,astolfi2010model}. The main idea is to use (rational) interpolation theory to approximate the transfer function of a system by another transfer function of lower degree.  Moment matching consists in imposing that the moments, \textit{i.e.}  the coefficients of the Laurent series expansion, of both transfer functions coincide at given points of the complex plane. Methods based on moment matching are numerically reliable and can be implemented efficiently using Krylov projectors~\cite[Chapter 11]{antoulas2005approximation}. Over the past two decades, moments of linear time-invariant systems have been characterized in terms of Sylvester equations~\cite{gallivan2004sylvester,gallivan2006model} and, under certain assumptions, in terms of steady-state responses~\cite{astolfi2010model}. These characterizations, in turn, have led to new model reduction methods for nonlinear systems~\cite{astolfi2010model}, for systems in explicit form~\cite{scarciotti2015modelexplicit}, for time-delay systems~\cite{scarciotti2016model}, and for systems with isolated singularities~\cite{padoan2017poles,padoan2017mrp,padoan2017eigenvalues,padoan2019isolated1}.   

However, a significant limitation of these methods is that the interpolation conditions imposed by moment matching are required to hold \emph{exactly}, which, for some purposes, is an unnecessarily stringent assumption. In practice, one can often tolerate a small error around each interpolation point and seek for the ``best'' model which minimizes these errors. Another limitation is the lack of error bounds, which precludes any \textit{a priori} guarantee on the quality of approximation.  

Least squares moment matching provides a particularly interesting solution to both issues~\cite{shoji1985model,aguirre1992least,aguirre1994partial,smith1995least,gugercin2006model}.  The main idea is to require that the interpolation conditions imposed by moment matching are satisfied only in a least squares sense. Model reduction methods based on least squares moment matching  
thus overcome the issues mentioned above by minimizing an optimization criterion, which directly yields \textit{a priori} error bounds and, under certain assumptions, guaranteed stability properties~\cite{gugercin2006model}. For linear systems, there is a vast literature on least squares moment matching~\cite{antoulas2020interpolatory,shoji1985model,aguirre1992least,aguirre1994partial,smith1995least,
gugercin2006model,gustavsen1999rational,
mayo2007framework,berljafa2017rkfit,nakatsukasa2018aaa}, with deep connections to Pad\'e approximation~\cite{aguirre1992least,aguirre1994partial} and Prony's method for filter design~\cite{gugercin2006model,parks1987digital} 
(see also~\cite{antoulas2005approximation,antoulas2020interpolatory}, and references therein, for further detail). However, model reduction by least squares moment matching does not have a nonlinear counterpart to date.

The goal of this work is to develop a unifying model reduction framework based on the notion of least squares moment matching for linear and nonlinear systems. The main ingredient of our approach is the formalism introduced in~\cite{astolfi2010model}, where 
moments of nonlinear systems have been defined and characterized using tools from output regulation theory~\cite{isidori1990output} (see also \cite[Chapter 8]{isidori1995nonlinear}).
The starting point of our analysis is a new time-domain characterization of least squares moment matching for linear systems, which relies on the solution of a constrained optimization problem involving a Sylvester equation. Building on the preliminary results presented in~\cite{padoan2021model}, models achieving least squares moment matching are shown to minimize an \textit{a priori} error bound expressed in terms of the worst case r.m.s. gain of an error system with respect to the family of signals produced by a given signal generator. Furthermore, a parameterized family of models achieving least squares moment matching is determined and shown to admit a natural system-theoretic interpretation in terms of a two-step model reduction process.

The main contribution of the paper is the development of analogous results for nonlinear systems. A nonlinear enhancement of the notion of least squares moment matching is \textit{defined} in terms of a constrained optimization problem involving an invariance equation. In close analogy with the linear case, models achieving least squares moment matching are then shown to minimize an \textit{a priori} error bound on the worst case r.m.s. gain of an error system with respect to the family of signals produced by a given signal generator. Finally, a parameterized family of models achieving least squares moment matching is presented and shown to possess an intuitive geometric interpretation, thus providing new insights on the linear theory.  

The proposed framework offers a unique (time-domain) perspective on the (frequency-domain) interpolation conditions imposed by least squares moment matching. 
In principle, our approach allows one to define least squares moment matching for systems which do not necessarily possess a representation in terms of transfer functions, including linear time-varying systems and hybrid systems, which we do not discuss for reason of space. Our approach also leads to a novel parameterization of models achieving least squares moment matching, which admits a natural system-theoretic interpretation in terms of a two-step model reduction process. This interpretation, in turn, bears significant practical implications and makes it possible to compute models achieving least squares moment matching using standard convex optimization solvers (such as CVX~\cite{cvx2014}).

The remainder of this work is organized as follows. 
Section~\ref{sec:linear} is devoted to the model reduction problem by least squares moment matching for linear systems. The problem is formulated and solved adopting a viewpoint which allows for a direct nonlinear counterpart to be developed. Section~\ref{sec:nonlinear} presents a nonlinear enhancement of the linear theory,  providing a unifying framework for least squares moment matching and new insights on the linear case. 
The theory is illustrated by means of worked-out examples in Section~\ref{sec:examples}. Section~\ref{sec:conclusion} concludes the paper with a summary and an outlook to future research directions.   

\textbf{Notation:} $\Zplus$, $\R$ and $\Co$ denote the set of non-negative integers, of real numbers, and of complex numbers, respectively. 
$\Cominus$ and $\Cozero$ denote the set of complex numbers with negative real part and with zero real part, respectively.
$\iota$ denotes the imaginary unit. 
$e_k$ denotes the vector  with the $k$-th entry equal to one and all other entries equal to zero. 
$I$ denotes the identity matrix.
$J_0$ denotes the matrix with ones on the superdiagonal and zeros elsewhere. 
$J_{s^{\star}}$ denotes the Jordan block associated with the eigenvalue ${s^{\star}\in\Co}$, \emph{i.e.} ${J_{s^{\star}} = s^{\star} I + J_0}$.
$\spectrum{A}$ denotes the spectrum of the matrix ${A \in \R^{n \times n}}$.
$\Image M$ and $\ker M$  denote the image and the kernel of the matrix ${M \in \R^{p \times m}}$, respectively. 
$M^{\transpose}$ and $M^{\dagger}$ denote the transpose and the Moore-Penrose inverse of the matrix ${M \in \R^{p \times m}}$, respectively. 
$\norm{\,\cdot\,}_2$ and $\norm{\,\cdot\,}_{2*}$ denote the Euclidean $2$-norm on $\R^{n}$ and the corresponding dual norm~\cite[p.637]{boyd2004convex}, respectively. Finally, $f^{(k)}(\cdot)$ denotes the derivative of order ${k\in\Zplus}$ of the function $f(\cdot)$, provided it exists, with $f^{(0)}(\cdot) = f(\cdot)$ by convention.

\section{Linear systems} \label{sec:linear}

\subsection{Problem formulation} \label{ssec:problem-formulation-linear}

Consider a continuous-time, single-input, single-output, linear, time-invariant system described by the equations
\beq \label{eq:system-linear}
\quad \dot{x} = Ax+Bu, \quad y=Cx,
\eeq
in which ${x(t)\in\R^n}$, ${u(t)\in\R}$, ${y(t)\in\R}$ and ${A\in\R^{n \times n}}$, ${B\in\R^{n \times 1}}$ and ${C\in\R^{1\times n}}$ are constant matrices, with transfer function defined as 
$${W(s)=C(sI-A)^{-1}B. }$$
For the notion of moment to make sense, we make the following standing assumption throughout the paper.
\begin{standing assumption}
The system~\eqref{eq:system-linear} is minimal, \textit{i.e.} controllable and observable.
\end{standing assumption}

\begin{definition}   \cite[p.345]{antoulas2005approximation}    \label{def:moment}
The \emph{moment of order ${k \in \Zplus}$} of system~\eqref{eq:system-linear} at ${s^{\star} \in \Co}$, with ${s^{\star} \not \in  \spectrum{A}}$, is defined as the complex number
\beq \nn
\eta_k(s^{\star}) = \frac{(-1)^k}{k!} W^{(k)}(s^{\star}) .
\eeq 
\end{definition}

\noindent 
Given distinct \emph{interpolation points} ${\{s_i\}_{i=1}^N}$, with 
${s_i \in \Co}$ and ${s_i \not \in  \spectrum{A}}$,
and the corresponding \emph{orders of interpolation}
${\{k_i\}_{i=1}^N}$, with ${k_i\in \Zplus}$,
model reduction by moment matching  consists in finding a system
\beq \label{eq:system-rom}
  \quad  \dot{\xi} = F\xi+Gv, \quad \psi = H\xi,
\eeq 
where ${\xi(t) \in\R^{r}}$, ${v(t)\in\R}$, ${\psi(t)\in\R}$ and ${F\in\R^{r \times r}}$, ${G\in\R^{r \times 1}}$ and ${H\in\R^{1\times r}}$ are constant matrices, the transfer function of which
\beq \nn
\hat{W}(s)=H(sI-F)^{-1}G 
\eeq
satisfies the \emph{interpolation conditions}
\beq  \label{eq:matching-condition-linear}
\quad \eta_{j}(s_i)= \hat{\eta}_{j}(s_i) ,\quad j \in \{0, \ldots, k_i\}, \ i \in \{1, \ldots, N\},
\eeq 
where $\eta_{j}(s_i) $ and $ \hat{\eta}_{j}(s_i)$ denote the moments of order $j$ of the systems~\eqref{eq:system-linear} and~\eqref{eq:system-rom} at $s_i$, respectively. The system~\eqref{eq:system-rom} is referred to as a \emph{model of system~\eqref{eq:system-linear}} and is said to \emph{achieve moment matching (at $\{ s_i\}^{N}_{i=1}$)} if the interpolation conditions~\eqref{eq:matching-condition-linear} hold~\cite[Chapter 11]{antoulas2005approximation}. Furthermore, if ${r < n}$, then~\eqref{eq:system-rom} is said to be a \emph{reduced order model of system~\eqref{eq:system-linear}}.

\begin{figure}[t!]
\centering
\includegraphics[width=\columnwidth,  height=5cm, keepaspectratio]{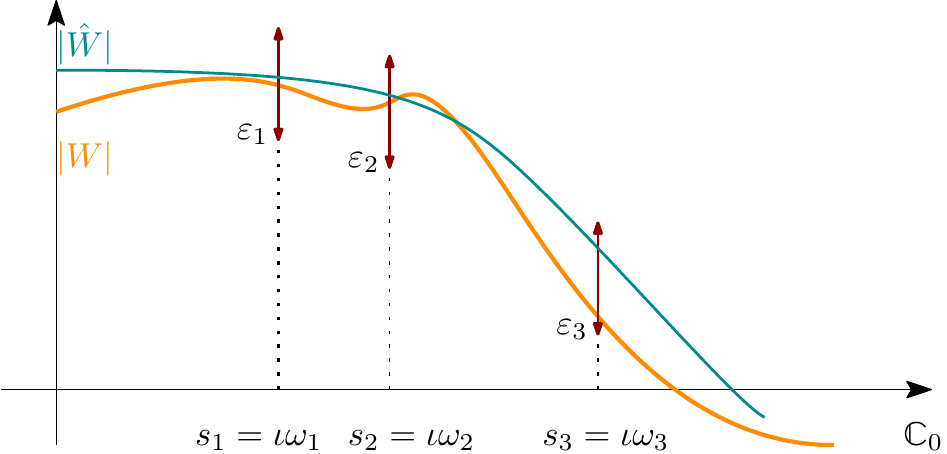}
\caption{The system~\eqref{eq:system-linear} must be approximated around every interpolation point $s_i =\iota\omega_i$ and a small error (depicted as $\varepsilon_i$ for the sake of illustration) can be tolerated. The order $r=1$ of the model is given and the original system is approximated at $N=3$ interpolation points.}
\label{fig:1}
\end{figure}

We are interested in the following problem. Suppose we wish to approximate system~\eqref{eq:system-linear} with the model~\eqref{eq:system-rom} 
by moment matching at a given set of interpolation points  $\{ s_i\}^{N}_{i=1}$ with corresponding orders of interpolation $\{ k_i\}^{N}_{i=1}$. Assume that the number of interpolation conditions 
$$\nu = \sum_{i=1}^{N} (k_i+1)$$
is much larger than the order of system~\eqref{eq:system-rom}, \emph{i.e.} $\nu \gg r$, and that  around every interpolation point $s_i$ a small error can be tolerated. The goal is to construct the ``best'' model of system~\eqref{eq:system-linear} of order $r$ (in a sense to be made precise in the sequel). Figure~\ref{fig:1} provides a pictorial representation of the problem, where a model of order $r=1$ needs to be constructed to approximate the original system at $N=3$ interpolation points.  

It is well-known that a model of order $r$ can match up to $2r$ moments~\cite[Chapter 11]{antoulas2005approximation}. The number of interpolation conditions is thus larger than the number of moments that can be matched if $\nu > 2r$. In this case, the interpolation conditions~\eqref{eq:matching-condition-linear} give rise to an overdetermined system of equations which can be solved in a least squares sense, leading directly to the \emph{model reduction problem by least squares moment matching} for linear systems.

\begin{problem}
Consider system~\eqref{eq:system-linear}. Let $\{ s_i\}^{N}_{i=1}$ be a set of interpolation points, with ${s_i \in \Co}$ and ${s_i \not \in  \spectrum{A}}$, and let $\{ k_i\}^{N}_{i=1}$ be the corresponding orders of interpolation, with ${k_i\in \Zplus}$. Let ${\nu = \sum_{i=1}^{N} (k_i+1)}$ and ${r\in \Zplus}$, with ${2r < \nu}$.  Find, if possible, a system~\eqref{eq:system-rom} of order $r$ which minimizes the index
\beq \label{eq:index}
\mathcal{J} = \sum_{i=1}^{N} \sum_{j=0}^{k_i} \left| \eta_{j}(s_i) - \hat{\eta}_{j}(s_i) \right|^2 .
\eeq
\end{problem}

\begin{figure}[h!]
\centering
\includegraphics[scale=0.8,keepaspectratio]{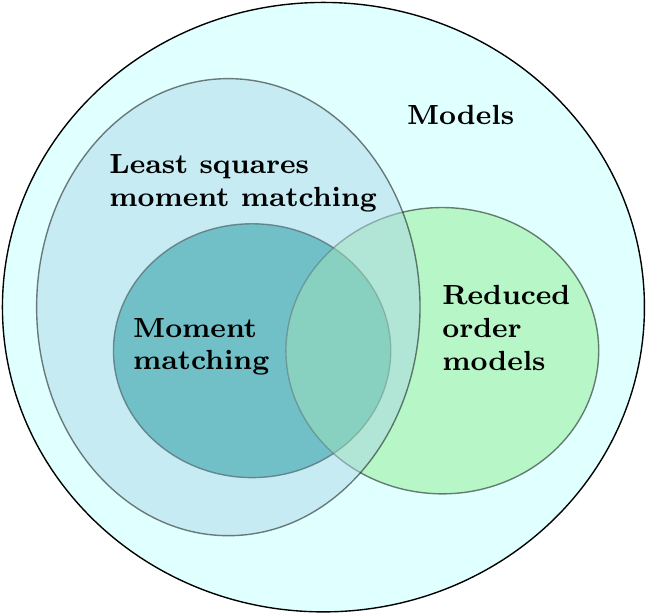}
\caption{Diagrammatic representation of the relationship among the notions of model, reduced order model, moment matching and least squares moment matching, respectively.}
\label{fig:2}
\end{figure}

For linear systems, the model reduction problem by least squares moment matching can be seen as a rational (Hermite) interpolation problem~\cite{davis1975interpolation} and, thus, it is generically well-posed (see, e.g.,~\cite{gutknecht1990sense} for further detail). The model~\eqref{eq:system-rom} is said to \emph{achieve least squares moment matching (at $\{ s_i\}^{N}_{i=1}$)} if it minimizes the index~\eqref{eq:index}. Clearly, any model achieving moment matching also achieves least squares moment matching, since the index~\eqref{eq:index} is minimized if the interpolation conditions~\eqref{eq:matching-condition-linear} hold. The set of (reduced order) models of system~\eqref{eq:system-linear} achieving moment matching is therefore a strict subset of the set of (reduced order) models of system~\eqref{eq:system-linear} achieving least squares moment matching. 
The Venn diagram in Figure~\ref{fig:2} summarizes the relationship among the notions of model, reduced order model, moment matching and least squares moment matching, respectively.

\subsection{Moment matching} \label{ssec:preliminaries}

Our analysis is relies on standard tools to achieve moment matching, including the use of Sylvester equations and Krylov projectors, which we recall below from~\cite{grimme1997krylov,gallivan2004sylvester,gallivan2006model,astolfi2010model} and~\cite[Chapter 11]{antoulas2005approximation} with minor variations.

\subsubsection{Moment matching via Sylvester equations}

\begin{lemma} \cite[Lemma 2]{astolfi2010model}  \label{lemma:astolfi-2}
Consider system~\eqref{eq:system-linear}.  Let ${k \in \Zplus}$ and ${s^{\star} \in \Co}$, with ${s^{\star} \not \in  \spectrum{A}}$. Then
$${
C\Pi \Psi 
=
\left[ \, \eta_0(s^{\star}) 	\; \eta_1(s^{\star}) \; \cdots \; \eta_k(s^{\star}) \, \right],
}$$
with ${\Psi\in \R^{(k+1) \times (k+1)}}$ a signature matrix\footnote{A signature matrix is a diagonal matrix with $\pm 1$ on the main diagonal~\cite[p.44]{bapat2010graphs}.} and ${\Pi   \in \R^{n\times (k+1) }}$ the unique solution of the Sylvester equation
$${
A\Pi  +B e_1^{\transpose} =\Pi   J_{s^{\star}} .
}$$
\end{lemma}

\begin{lemma} \cite[Lemmas 3 and 4]{astolfi2010model}   \label{lemma:astolfi-3}
Consider system~\eqref{eq:system-linear}. 
Assume ${S\in \R^{\nu \times \nu}}$ is a non-derogatory\footnote{A matrix is non-derogatory if its characteristic polynomials and its minimal polynomial coincide~\cite[p.178]{horn1994matrix}.} matrix with characteristic polynomial 
\beq \label{eq:characteristic-polynomial}
\chi_S(s)=\prod_{i=1}^{N} (s-s_i)^{k_i+1}
\eeq
and ${L\in\R^{1\times \nu}}$ is such that the pair $(S,L)$ is observable.
Then the moments $\eta_0(s_1),$ $\dots,$ $\eta_{k_1}(s_1)$, $\dots,$ $\eta_0(s_N),$ $\dots,$ $\eta_{k_N}(s_N)$ are in one-to-one correspondence\footnote{The terminology is borrowed from~\cite{astolfi2010model}, where one-to-one correspondence between two objects means that one uniquely determines the other and \emph{vice versa}.} with the matrix ${C\Pi}$, where  ${\Pi \in \R^{n\times \nu}}$ is the unique solution of the Sylvester equation
\beq \label{eq:Sylvester-equation-astolfi}
A\Pi +BL=\Pi S.
\eeq
\end{lemma}%
\noindent
Lemma~\ref{lemma:astolfi-3} establishes that the moments of system~\eqref{eq:system-linear} can be \emph{equivalently} characterized in terms of the product of the output matrix of the system and the solution of the Sylvester equation~\eqref{eq:Sylvester-equation-astolfi}.  This, in turn, is instrumental to provide a time-domain characterization of the moments of system~\eqref{eq:system-linear} in terms of the \emph{steady-state output response}\footnote{See~\cite[Chapter 8]{isidori1995nonlinear} for the notion of steady-state response.} of the interconnection of system~\eqref{eq:system-linear} with a signal generator described by the equations
\beq \label{eq:system-signal-generator}
 \dot{\omega} = S\omega, \quad  \theta = L\omega,  
\eeq 
with ${\omega(t) \in \R^\nu}$ and ${\theta(t) \in \R}$, which satisfies the following assumptions.

\begin{assumption} \label{ass:signal-generator-linear}
The matrix ${S\in \R^{\nu \times \nu}}$ is non-derogatory and has characteristic polynomial~\eqref{eq:characteristic-polynomial}. The matrix ${L\in\R^{1\times \nu}}$ is such that the pair $(S,L)$ is observable.
\end{assumption}

\begin{assumption} \label{ass:excitability-linear}
The vector ${\omega(0) \in\R^{\nu}}$ is such that the pair $(S,\omega(0))$ is controllable.
\end{assumption}

\begin{theorem} \cite[Theorem 1]{astolfi2010model}   \label{thm:astolfi-1}
Consider system~\eqref{eq:system-linear} and the signal generator~\eqref{eq:system-signal-generator}.
Suppose Assumptions~\ref{ass:signal-generator-linear} and~\ref{ass:excitability-linear} hold.
Assume ${\spectrum{A} \subset \Cominus}$ and ${\spectrum{S} \subset \Cozero}$. Then the moments $\eta_0(s_1),$ $\dots,$ $\eta_{k_1}(s_1)$, $\dots,$ $\eta_0(s_N),$ $\dots,$ $\eta_{k_N}(s_N)$ are in one-to-one correspondence with  the steady-state output response of the interconnected system~\eqref{eq:system-linear}-\eqref{eq:system-signal-generator}, with ${u =  \theta}$.
\end{theorem}

\noindent
Theorem~\ref{thm:astolfi-1} motivates the following notion of model achieving moment matching.

\begin{definition}  \cite[p.4]{astolfi2010model} \label{def:reduced-order-model-linear}
The system~\eqref{eq:system-rom} is a \emph{model of system~\eqref{eq:system-linear} at $(S,L)$}, with $S \in \R^{\nu\times \nu}$ a non-derogatory matrix such that ${\spectrum{A} \cap \spectrum{S} = \emptyset}$, if ${\spectrum{F} 
\cap \spectrum{S} = \emptyset}$
and 
\beq \label{eq:model-condition-2}
C\Pi = HP,
\eeq
where ${\Pi \in \R^{ n \times \nu}}$ is the unique solution of the Sylvester equation
\eqref{eq:Sylvester-equation-astolfi}, ${L \in \R^{1 \times \nu}}$ such that the pair $(S,L)$ is observable, and ${P \in \R^{ r \times \nu}}$ is the unique solution of the Sylvester equation
\beq \label{eq:Sylvester-equation-astolfi-model} 
FP +GL= P S.
\eeq
In this case, system~\eqref{eq:system-rom} is said to \emph{match the moment of system~\eqref{eq:system-linear} (or to achieve moment matching) at $(S,L)$}. Furthermore, system~\eqref{eq:system-rom} is a \emph{reduced order model of system~\eqref{eq:system-linear} at $(S,L)$} if ${r<n}$. 
\end{definition}

\begin{remark}
The solution ${P\in\R^{r \times \nu}}$ of the Sylvester equation~\eqref{eq:Sylvester-equation-astolfi-model}
is necessarily such that the subspace $\ker P$ is $(S,L)$ conditioned invariant\footnote{See~\cite[p.199]{basile1992controlled} for the definition of $(A,C)$ conditioned invariant subspace.}.
This follows immediately from~\cite[Property 4.1.4]{basile1992controlled}. 
\end{remark}

\noindent
A family of reduced order models achieving moment matching for system~\eqref{eq:system-linear} has been determined in~\cite{astolfi2010model} by selecting $r=\nu$ and ${P=I}$ in~\eqref{eq:model-condition-2} and~\eqref{eq:Sylvester-equation-astolfi-model}, respectively. As a matter of fact, this yields a family of reduced order models of system~\eqref{eq:system-linear} at $(S,L)$ described by equations~\eqref{eq:system-rom}, with
\beq \label{eq:family-linear-astolfi}
F=S-\Delta L, 
\quad G=\Delta, 
\quad  H = C \Pi,
\eeq
in which ${\Delta \in \R^{r \times 1}}$ is such that ${\spectrum{S-\Delta L} \cap \spectrum{S}  = \emptyset}$. The vector ${\Delta }$ is a \emph{parameter} of the family of reduced order models~\eqref{eq:family-linear-astolfi} which can be used to assign prescribed properties to the reduced order model, including stability, passivity, and a given $L_2$-gain~\cite{astolfi2010model}.

\subsubsection{Moment matching via Krylov methods} 

Krylov methods provide a numerically efficient way to construct models achieving moment matching~\cite[Chapter 11]{antoulas2005approximation}. Given system~\eqref{eq:system-linear}, the main idea is to define a model~\eqref{eq:system-rom}, by projection, as
\beq \label{eq:family-linear-Krylov}
F = P A Q,
\quad G=P B, 
\quad  H = C Q ,
\eeq
where ${P \in \R^{r \times n}}$ and ${Q \in \R^{n \times r}}$ are (biorthogonal Petrov-Galerkin) projectors such that $PQ=I$ and such that the interpolation conditions~\eqref{eq:matching-condition-linear} hold~\cite[Chapter 11]{antoulas2005approximation}.  The key tool to construct the projectors $P$ and $Q$ is the notion of \emph{Krylov subspace}, which for a matrix ${M \in \Co^{n \times n}}$, a vector ${v\in\Co^n}$, a point ${s^{\star} \in \Co \cup \{ \infty\}}$, and an index $j\in\Zplus$, with $j>0$, 
is defined, for ${s^{\star} = \infty}$, as 
\beq \nn 
\mathcal{K}_j(M,v; s^{\star}) =     \Span \{ \, v , \,  \ldots , \, M^{j-1} v \,\}      , 
\eeq 
and,  for ${s^{\star} \not = \infty}$, as
\beq \nn
\mathcal{K}_j(M,v; s^{\star}) =      \Span \{  (s^{\star}I-M)^{-1} v  ,\, \ldots ,\, (s^{\star}I-M)^{-j} v \}     .
\eeq
Models achieving moment matching can be constructed, for example, by ensuring that the kernel of the projector $P$ is contained in the intersection of certain Krylov subspaces, as detailed by the following statement.

\begin{theorem}~\cite[Theorem 3.1]{grimme1997krylov} \label{thm:Krylov}
Consider system~\eqref{eq:system-linear}.  Let $\{ s_i\}^{N}_{i=1}$ be a set of interpolation points, with ${s_i \in \Co}$ and ${s_i \not \in  \spectrum{A}}$, and let $\{ k_i\}^{N}_{i=1}$ be the corresponding orders of interpolation, with ${k_i\in \Zplus}$. Define the model~\eqref{eq:system-rom} as in~\eqref{eq:family-linear-Krylov}, with ${P \in \R^{r \times n}}$ and ${Q \in \R^{n \times r}}$ such that $PQ=I$ and such that
$${\ker P \subset \bigcap_{i =1}^{N} \mathcal{K}_{k_i}(A^{\transpose},C^{\transpose}; s_i)}.$$ 
Then the model~\eqref{eq:system-rom} achieves moment matching at $\{ s_i\}^{N}_{i=1}$ .
\end{theorem}

\noindent
For more detail on model reduction via Krylov methods, the reader is referred to~\cite{grimme1997krylov,gallivan2004sylvester,gallivan2006model} and~\cite[Chapter 11]{antoulas2005approximation}.

\subsection{Least squares moment matching}

We begin our analysis by establishing that models achieving least squares moment matching can be \textit{equivalently} characterized in terms of the solutions of a constrained optimization problem of the form
\beq \label{eq:optimization-problem}
\begin{array}{ll}
    \mbox{minimize}      & \norm{(C \Pi - H P)T }_{2*}^2 \\
    \mbox{subject to}    & FP +GL= P S , \\
    						& \spectrum{F} \cap \spectrum{S} = \emptyset,
\end{array}
\eeq
for some non-singular matrix ${T \in \R^{\nu\times \nu}}$, where ${F\in\R^{r \times r}}$, ${G\in\R^{r \times 1}}$, ${H\in\R^{1\times r}}$ and ${P\in\R^{r\times \nu}}$  are the optimization variables, while system~\eqref{eq:system-linear} and the signal generator~\eqref{eq:system-signal-generator} 
(and, thus, the solution ${\Pi \in \R^{ n \times \nu}}$ of the Sylvester equation~\eqref{eq:Sylvester-equation-astolfi}) 
are problem data. To this end, we first introduce a basic assumption and prove a preliminary lemma, which allows us to rewrite the index~\eqref{eq:index} in terms of the solutions of the Sylvester equations~\eqref{eq:Sylvester-equation-astolfi} and~\eqref{eq:Sylvester-equation-astolfi-model}.
 
\begin{assumption} \label{ass:T-linear}
The matrix ${T \in \R^{\nu\times \nu}}$ is non-singular and such that
\beq \label{eq:thm-index-proof-00}
ST = TJ, \quad L T = \Lambda ,
\eeq
with ${J = \diag(J_{s_1}, \ldots, J_{s_N})}$ and $\Lambda = [\,e_1^{\transpose} \, \cdots \, e_1^{\transpose} \,]$. 
\end{assumption}

\begin{lemma} \label{lem:index}
Consider system~\eqref{eq:system-linear}, the model~\eqref{eq:system-rom} and the signal generator~\eqref{eq:system-signal-generator}.  Suppose Assumptions~\ref{ass:signal-generator-linear} and~\ref{ass:T-linear} hold. Assume ${\spectrum{A}  \cap \spectrum{S} = \emptyset}$ and ${ \spectrum{F} \cap \spectrum{S} = \emptyset}$. Then 
\beq \label{eq:thm-index-inequality}
\mathcal{J} =  \norm{(C \Pi - H P)T }_{2*}^2 , 
\eeq
where ${\Pi \in\R^{n \times \nu}}$ and ${P\in\R^{r \times \nu}}$ are the (unique) solutions of the Sylvester equations~\eqref{eq:Sylvester-equation-astolfi} and~\eqref{eq:Sylvester-equation-astolfi-model}, respectively.
\end{lemma}

\begin{proof}
The assumptions $\spectrum{A}  \cap \spectrum{S} = \emptyset$ and $ \spectrum{F} \cap \spectrum{S} = \emptyset$ directly imply the existence and uniqueness of the solutions of the Sylvester equations~\eqref{eq:Sylvester-equation-astolfi} and~\eqref{eq:Sylvester-equation-astolfi-model}, respectively. Furthermore, in view of Assumptions~\ref{ass:signal-generator-linear} and~\ref{ass:T-linear}, \eqref{eq:Sylvester-equation-astolfi} and~\eqref{eq:thm-index-proof-00} together imply
\beq \nn
A\Pi +B\Lambda T^{-1}=\Pi TJT^{-1},
\eeq
or, equivalently,
\beq \nn
A\Pi T +B\Lambda =\Pi TJ ,
\eeq
where ${\Pi \in \R^{n \times \nu}}$ is the (unique) solution of the Sylvester equation~\eqref{eq:Sylvester-equation-astolfi}. Then, setting  ${\widetilde \Pi = \Pi T}$ and appealing to Lemma~\ref{lemma:astolfi-2}, one obtains
\beq \label{eq:thm-index-proof-01}
C\widetilde \Pi \Psi \!=\! \left[ \, \eta_0(s_1) 	\,  \cdots \, \eta_{k_1}(s_1) \, \cdots \, \eta_0(s_N) 	\,  \cdots \, \eta_{k_N}(s_N) \, \right]  \!,
\eeq
in which ${\Psi \in \R^{\nu \times \nu}}$ is a signature matrix. A similar reasoning applies to the model~\eqref{eq:system-rom}. Thus, setting ${\widetilde P = P T}$ and appealing to Lemma~\ref{lemma:astolfi-2}, one obtains 
\beq \label{eq:thm-index-proof-02}
H\widetilde P  \Psi \!=\! \left[ \, \hat\eta_0(s_1) 	\, \cdots \, \hat\eta_{k_1}(s_1) \, \cdots \, \hat\eta_0(s_N) 	\,  \cdots \, \hat\eta_{k_N}(s_N) \, \right] \!,
\eeq 
where ${P\in \R^{r \times \nu}}$ is the (unique) solution of the Sylvester equation~\eqref{eq:Sylvester-equation-astolfi-model}.
Then 
\beq \nn
\scalebox{0.95}{$
\begin{array}{rcl}
\norm{(C \Pi - H P)T }_{2*}^2 
&=& \norm{C\widetilde \Pi  - H \widetilde P }_{2*}^2  \\
&=& \norm{C\widetilde \Pi \Psi - H \widetilde P \Psi }_{2*}^2  \\
&\stackrel{\eqref{eq:thm-index-proof-01},\eqref{eq:thm-index-proof-02}}{=}& 
\displaystyle
\sum_{i=1}^{N} \sum_{j=0}^{k_i} \left| \eta_{j}(s_i) - \hat{\eta}_{j}(s_i) \right|^2 ,
\end{array}
$}
\eeq
where the last identity holds since $\Psi$ is a signature matrix and since the norm $\norm{\,\cdot\,}_{2*}$ is unitarily invariant.
\end{proof}%

\begin{remark} \label{rem:T-unitary}
Assumption~\ref{ass:signal-generator-linear} requires that the pair $(S,L)$ is observable. This guarantees the existence of a matrix ${T \in \R^{\nu\times \nu}}$ such that~\eqref{eq:thm-index-proof-00} holds. However, the proof of Lemma~\ref{lem:index} shows that the matrix $T$ acts as a \textit{weight} on the norm of the covector ${C \Pi - H P}$. This, in turn, suggests that the least squares approximation can be modulated in regions of interest by relaxing Assumptions~\ref{ass:T-linear} and by requiring instead that the matrix ${T}$  
is merely non-singular. For reasons of space, we do not elaborate further on
the role of weights, but defer a full investigation of this issue to future work.
\end{remark}

We are now ready to show that, under certain assumptions, the interpolation constraints imposed by least squares moment matching can be \textit{equivalently} characterized in terms of the solutions of the optimization problem~\eqref{eq:optimization-problem}.

\begin{theorem} \label{thm:moments-optimization-linear}
Consider system~\eqref{eq:system-linear}, the model~\eqref{eq:system-rom} and the signal generator~\eqref{eq:system-signal-generator}.  Suppose Assumptions~\ref{ass:signal-generator-linear} and~\ref{ass:T-linear} hold. Assume ${\spectrum{A}  \cap \spectrum{S} = \emptyset}$. Then the model~\eqref{eq:system-rom} achieves least squares moment matching at $\spectrum{S}$ if and only if there exists a matrix ${P\in\R^{r\times \nu}}$ such that $(F,G,H,P)$ is a solution of the optimization problem~\eqref{eq:optimization-problem}. Furthermore, the model~\eqref{eq:system-rom} is controllable if and only if the matrix $P$ is full rank.
\end{theorem}

\begin{proof}
Under the stated assumptions, the moments of the model~\eqref{eq:system-rom} at $\spectrum{S}$ are well-defined if and only if ${\spectrum{F} \cap \spectrum{S} = \emptyset}$, which directly implies the existence and uniqueness of a solution $P$ of the Sylvester equation~\eqref{eq:Sylvester-equation-astolfi-model}. The first claim thus follows immediately from Lemma~\ref{lem:index} and from the definition of model achieving least squares moment matching. 

We now show that the model~\eqref{eq:system-rom} is controllable if and only if the matrix $P$ is full rank.
To this end, assume, by contraposition, that $(F,G)$ is not controllable. The Popov-Belevitch-Hautus criterion implies that there exists ${\lambda \in \Co}$ and ${w\in\Co^r}$, with ${w\not  = 0}$, such that $w^{\transpose} [\, \lambda I - F \, | \, G \,] = 0$ or, equivalently, 
\beq \label{eq:moments-optimization-linear-proof}
\lambda w^{\transpose} = w^{\transpose}  F, \quad w^{\transpose}G = 0.  
\eeq
Then 
\beq \nn
w^{\transpose} P S 
 \stackrel{\eqref{eq:Sylvester-equation-astolfi-model}}{=}   w^{\transpose} (FP+GL) 
 =  w^{\transpose}FP+w^{\transpose}GL 
 \stackrel{\eqref{eq:moments-optimization-linear-proof}}{=}  \lambda w^{\transpose}P,
\eeq
which, in view of \cite[Theorem 1]{de1981controllability}, yields 
$${\rank P < r,}$$
since, by assumption, the pair $(S,L)$ is observable and ${\spectrum{F} \cap \spectrum{S} = \emptyset}$. 
The converse implication is easily obtained by reversing the above arguments.
\end{proof}

Theorem~\ref{thm:moments-optimization-linear} bears a number of interesting consequences. The connection
between least squares moment matching and the solutions of the optimization problem~\eqref{eq:optimization-problem} allows one to reformulate the interpolation constraints~\eqref{eq:matching-condition-linear} imposed by least squares moment matching in terms of the solutions of the Sylvester equations~\eqref{eq:Sylvester-equation-astolfi} and~\eqref{eq:Sylvester-equation-astolfi-model}. This formulation, in turn, is instrumental to \textit{define} a notion of least squares moment matching for nonlinear systems which do not necessarily possess a representation in terms of a transfer function, as discussed in detail in Section~\ref{sec:nonlinear}.

Another interesting consequence of Theorem~\ref{thm:moments-optimization-linear} concerns the solution of the least squares model reduction problem. In principle, the problem can be solved in two steps. The first step consists in solving the Sylvester equation~\eqref{eq:Sylvester-equation-astolfi} in the unknown matrix $\Pi$. The second step consists in computing a solution of the optimization problem~\eqref{eq:optimization-problem} in the unknown matrices $F,$ $G,$ $H,$ and $P$, while possibly enforcing specific properties to the model sought. 
This computation can be carried out in the spirit of~\cite{astolfi2010model} by regarding the optimization variable $P$ as a given \emph{parameter}. This leads to a family of models achieving least squares moment matching, which admits a natural system-theoretic interpretation and which directly relates the parameters to system-theoretic properties of the model. We elaborate more on these ideas in Sections~\ref{ssec:LS-MR-MM-linear},~\ref{ssec:system-theoretic-interpretation} and~\ref{ssec:properties}.

Finally, Theorem~\ref{thm:moments-optimization-linear} allows one to infer an \textit{a priori} error bound expressed in terms of the worst case r.m.s. gain of an error system with respect to the family of signals produced by a given signal generator, as detailed in the next subsection.

\subsubsection{An \textit{a priori} error bound}

Least squares moment matching can be given a simple interpretation in terms of the steady-state behavior of the error system 
\beq \label{eq:system-error-linear}
 \dot{x} = Ax+B u,  \quad 
 \dot{\xi} = F\xi+G u,   \quad 
 e=Cx-H\xi,  
\eeq 
in which  ${x(t)\in\R^n}$, ${\xi(t)\in\R^r}$, ${u(t)\in\R}$, and ${e(t)\in\R}$. In particular, if all solutions of the signal generator~\eqref{eq:system-signal-generator} are periodic and if the steady-state output response $e_{ss}$ of the  interconnected system~\eqref{eq:system-signal-generator}-\eqref{eq:system-error-linear}, with $u=\theta$, is well-defined, then achieving least squares moment matching corresponds to minimizing an upper bound of the \textit{worst case r.m.s. gain} of the error system~\eqref{eq:system-error-linear}  with respect  to the family of signals produced by the signal generator~\eqref{eq:system-signal-generator}, defined as~\cite[p.98]{boyd1991linear}
\beq \label{eq:gain-rms-linear}
\gamma_{rms} =  \sup_{\omega(\cdot) \in \mathcal{W}} \frac{\norm{e_{ss}}_{rms}}{\norm{\omega}_{rms}}
\eeq 
where $\norm{\omega}_{rms}$ is the \textit{r.m.s. value} of the signal ${\omega(t) \in \R^\nu}$, defined as~\cite[p.86]{boyd1991linear}
\beq \label{eq:rms}
\norm{\omega}_{rms} = \left(\lim_{\tau \to \infty} \frac{1}{\tau} \int_{0}^{\tau} \norm{\omega(t)}_2^2 dt \right)^{1/2} ,
\eeq
provided the limit exists, while $\mathcal{\mathcal{W}}$ is the family of signals produced by~\eqref{eq:system-signal-generator} with non-zero r.m.s. value.

\begin{theorem} \label{thm:rms-linear}
Consider system~\eqref{eq:system-linear}, the model~\eqref{eq:system-rom} and the signal generator~\eqref{eq:system-signal-generator}. Suppose Assumptions~\ref{ass:signal-generator-linear}-\ref{ass:T-linear} hold. 
Assume ${\spectrum{A} \cup \spectrum{F} \subset \Cominus}$ and
${S + S^{\transpose} = 0}$. Then the following statements hold.
\begin{itemize}
\item[(i)] The steady-state output response of the interconnected system~\eqref{eq:system-signal-generator}-\eqref{eq:system-error-linear}, with ${u=\theta}$, is well-defined and uniquely determined by the moments of the error system~\eqref{eq:system-error-linear}  at $\spectrum{S}$.
\item[(ii)] The worst case r.m.s. gain of the error system~\eqref{eq:system-error-linear} 
 with respect  to the family of signals produced by the signal generator~\eqref{eq:system-signal-generator}  is well-defined and such that
\beq \label{eq:gain-bound-linear}
\gamma_{rms}  \le \norm{C \Pi - H P}_{2*} ,
\eeq
with ${\Pi \in\R^{n \times \nu}}$ and ${P\in\R^{r \times \nu}}$ the (unique) solutions of the Sylvester equations~\eqref{eq:Sylvester-equation-astolfi} and~\eqref{eq:Sylvester-equation-astolfi-model}, respectively. 
\item[(iii)]  The error bound~\eqref{eq:gain-bound-linear} is minimized if the model~\eqref{eq:system-rom} achieves least squares moment matching at~$\spectrum{S}$ and if the matrix ${T\in\R^{\nu \times \nu}}$ is orthogonal\footnote{A matrix ${T \in \R^{n \times n}}$ is orthogonal if ${TT^{\transpose} =  I}$ \cite[p.84]{horn1994matrix}.}.
\end{itemize}
\end{theorem}

\begin{proof}
The proof relies on arguments which can be naturally transposed to nonlinear systems. A simpler, alternative
proof can be obtained via Laplace transforms.

(i). The dynamics of the interconnected system~\eqref{eq:system-signal-generator}-\eqref{eq:system-error-linear}, with $u=\theta$, is governed by the equations
\beq \label{eq:system-error-interconnected-linear} 
\dot{\omega} \!=\! S\omega, \
\dot{x} \!=\! Ax+B L\omega,  \
 \dot{\xi} \!=\! F\xi+G L\omega,   \
 e\!=\!Cx-H\xi .
\eeq 
By assumption, ${\spectrum{A} \cup \spectrum{F} \subset \Cominus}$. Moreover, ${\spectrum{S} \subset \Cozero}$, 
since ${S + S^{\transpose} = 0}$. By the center manifold theorem~\cite[p.4]{carr1981application}, the interconnected system~\eqref{eq:system-error-interconnected-linear} has therefore a well-defined center manifold
\beq \nn
\mathcal{M}_c = \left\{ \, (x,\xi,\omega) \,:\,  x = \Pi\omega, \, \xi =  P\omega \,   \right\} ,
\eeq
where ${\Pi \in\R^{n \times \nu}}$ and ${P\in\R^{r \times \nu}}$ are the (unique) solutions of the Sylvester equations~\eqref{eq:Sylvester-equation-astolfi} and~\eqref{eq:Sylvester-equation-astolfi-model}, respectively. 
Furthermore, the manifold $\mathcal{M}_c$ is exponentially attractive, since
\beq\nn
\scalebox{0.95}{$
\renewcommand{\arraystretch}{1}
\setlength{\arraycolsep}{1.5 pt}
\dot{\overbrace{\bma \begin{array}{c} x- \Pi\omega \\ \xi - P \omega \end{array} \ema}} 
\!\stackrel{\eqref{eq:system-error-interconnected-linear}}{=} \!
\bma \begin{array}{c} Ax + BL\omega - \Pi S\omega \\ F\xi + GL\omega - P S\omega \end{array} \ema 
\!\stackrel{\eqref{eq:Sylvester-equation-astolfi},\eqref{eq:Sylvester-equation-astolfi-model}}{=}\!
\bma \begin{array}{c} A(x- \Pi\omega) \\ F(\xi - P \omega) \end{array} \ema \!.$}
\eeq
The output of the interconnected system~\eqref{eq:system-error-interconnected-linear} can be thus written as
${e =  e_{ss} + e_{d}},$
in which
\beq \label{eq:error-steady-state-linear}
e_{ss} = (C\Pi - HP)\omega
\eeq
is the well-defined steady-state output response of the error system~\eqref{eq:system-error-linear} 
and $e_{d}$ is an exponentially decaying signal. Furthermore, under the stated assumptions, Theorem~\ref{thm:astolfi-1} implies that the steady-state output response $e_{ss}$ is uniquely determined by the moments of the error system~\eqref{eq:system-error-linear} at $\spectrum{S}$.

(ii). The worst case r.m.s. gain of the error system~\eqref{eq:system-error-linear} is  well-defined: by assumption, the matrix $S$ is skew-symmetric, which implies that the solutions of the signal generator~\eqref{eq:system-signal-generator} are periodic, and the steady-state output response of the interconnected system~\eqref{eq:system-error-interconnected-linear} is also periodic (with the same period) by~\eqref{eq:error-steady-state-linear}.   To obtain the error bound~\eqref{eq:gain-bound-linear}, note that
\begin{align}
\norm{e_{ss}(t)}_{2}  
&\stackrel{\eqref{eq:error-steady-state-linear}}{=} \!
\norm{(C\Pi - HP)\omega(t)}_{2} \nn \\
&\,\le
\norm{C\Pi - HP}_{2*}
\norm{\omega(t)}_{2} , \label{eq:coordinates-change-linear1}
\end{align}
where the last inequality follows from the Cauchy-Schwartz inequality. 
Then
\beq \nn
\gamma_{rms} 
\!\stackrel{\eqref{eq:gain-rms-linear}}{=} \! 
\sup_{\omega(\cdot) \in \mathcal{W}} \frac{\norm{e_{ss}}_{rms}}{\norm{\omega}_{rms}}
\!\stackrel{\eqref{eq:coordinates-change-linear1}}{\le} 
\norm{C\Pi - HP}_{2*} ,
\eeq
which proves the error bound~\eqref{eq:gain-bound-linear}. 

(iii) If the matrix $T$ is orthogonal, then 
\beq \nn
\norm{C\Pi - HP}_{2*} = \norm{(C\Pi - HP)T}_{2*},
\eeq
since the norm $\norm{\,\cdot\,}_{2*}$ is unitarily invariant. The claim is thus a direct consequence of item (ii) and Theorem~\ref{thm:moments-optimization-linear}.
\end{proof}

Theorem~\ref{thm:rms-linear} formalizes the idea that, under certain assumptions, model reduction by least squares  moment matching corresponds to minimizing an upper bound of the worst case r.m.s. gain of the error system~\eqref{eq:system-error-linear}  with respect  to the family of signals produced by the signal generator~\eqref{eq:system-signal-generator}. This intuitive steady-state interpretation of least squares moment matching admits a nonlinear analogue, as discussed in detail  in Section~\ref{sec:nonlinear}.

\begin{remark}
The existence of an orthogonal matrix ${T \in \R^{\nu\times \nu}}$ which brings the matrix $S$ in (real) Jordan form
is guaranteed if the matrix $S$ is normal\footnote{A matrix ${N \in \R^{n \times n}}$ is normal if ${NN^{\transpose} =  N^{\transpose}N}$~\cite[p.131]{horn1994matrix}.} (or, in particular, skew-symmetric).
\end{remark}

\subsection{Model reduction by least squares moment matching}  \label{ssec:LS-MR-MM-linear}

The results of the previous section can be used to obtain a family of models achieving least squares moment matching by regarding the optimization variable $P$ as a given \emph{parameter}.
The family of models in question is described by the equations~\eqref{eq:system-rom}, with 
\beq \label{eq:family-linear}
F=P(S-\Delta L)Q,
\quad G=P \Delta , 
\quad  H = C \Pi Q,
\eeq
in which ${S\in \R^{\nu \times \nu}}$ is a non-derogatory matrix with characteristic polynomial~\eqref{eq:characteristic-polynomial} such that ${\spectrum{A}  \cap \spectrum{S} = \emptyset}$, 
${L\in\R^{1\times \nu}}$ is such that the pair $(S,L)$ is observable,
${\Pi \in \R^{n \times \nu}}$ is the (unique) solution of the Sylvester equation~\eqref{eq:Sylvester-equation-astolfi},
while ${P \in \R^{r \times \nu}}$, ${\Delta \in\R^{\nu \times 1}}$ and ${Q \in\R^{\nu \times r}}$  are such that
\begin{itemize}
\item[(A$_P$)]  the matrix $P$ is full rank and such that the subspace ${\ker P}$ is  ${(S,L)}$ conditioned invariant,
\item[(A$_Q$)]  the matrix $Q$ is a right inverse
of $P$, \textit{i.e.} ${PQ=I}$, defined as  
\beq \label{eq:Q}
{Q = TT^{\transpose} P^{\transpose} (PTT^{\transpose} P^{\transpose})^{-1}},   
\eeq
with ${T\in\R^{\nu \times \nu}}$ any matrix such that Assumption~\ref{ass:T-linear} holds,
\item[(A$_\Delta$)] the vector $\Delta$ is such that the subspace~$\ker P$ is ${(S-\Delta L)}$-invariant and such that 
\beq
{ \spectrum{P(S-\Delta L)Q}  \cap \spectrum{S} = \emptyset},
\eeq
\end{itemize}
in which case $P$, $\Delta$ and $Q$ are said to be \emph{admissible} for the parameterization~\eqref{eq:family-linear}.

The family of models~\eqref{eq:system-rom}-\eqref{eq:family-linear} provides a solution to the model reduction problem by least squares moment matching, as detailed by the following statement.    

\begin{theorem} \label{thm:1}
Consider system~\eqref{eq:system-linear} and the family of models~\eqref{eq:system-rom}-\eqref{eq:family-linear}. 
Let ${P\in\R^{r\times \nu}}$,  ${Q\in\R^{\nu\times r}}$ and ${\Delta\in\R^{\nu \times 1}}$ 
be admissible for the parameterization~\eqref{eq:family-linear}. 
Then the model~\eqref{eq:system-rom}-\eqref{eq:family-linear}
achieves least squares moment matching at $\spectrum{S}$.
\end{theorem}

\begin{proof}
By Theorem~\ref{thm:moments-optimization-linear}, it suffices to establish that $(F,G,H,P)$, with $F$, $G$ and $H$ defined as in~\eqref{eq:family-linear}, is a solution of the optimization problem~\eqref{eq:optimization-problem} for some matrix ${P \in \R^{r\times \nu}}$ and some  matrix 
 ${T \in \R^{\nu\times \nu}}$ such that Assumption~\ref{ass:T-linear} holds. 
Equivalently, it suffices to show that all admissible parameters yield a solution of the optimization problem
\beq \label{eq:optimization-problem-proof}
\begin{array}{ll}
    \mbox{minimize}      & \norm{C \Pi(I-QP)T}_{2*}^2 \\
    \mbox{subject to}    & P(S-\Delta L)(I-QP) = 0, \\
    						&  \spectrum{P(S-\Delta L)Q} \cap  \spectrum{S}= \emptyset,
\end{array}
\eeq
with optimization variables $P$, $Q$ and $\Delta$,  since~\eqref{eq:optimization-problem-proof} can be obtained, after simple algebraic manipulations, by direct substitution of~\eqref{eq:family-linear} into~\eqref{eq:optimization-problem}.

Let ${P}$, ${\Delta}$ and ${Q}$ be admissible for the parameterization~\eqref{eq:family-linear} and let ${T}$ be such that Assumption~\ref{ass:T-linear} holds. Note that (A$_\Delta$) trivially implies the condition ${\spectrum{P(S-\Delta L)Q} \cap  \spectrum{S} = \emptyset}$. Furthermore, (A$_P$) and (A$_\Delta$) imply that the subspace ${\ker P}$ is $(S-\Delta L)$-invariant, \textit{i.e.} 
\beq \label{eq:lemma1-linear-1}
P(S-\Delta L)v =0,   \quad \forall \,  v \in \ker P.  
\eeq
Let ${\bar{v}_i  = (I - QP)e_i}$  for ${i \in \{1, \ldots, \nu\}}$. Note that ${\bar{v}_i \in \ker P}$, since 
$${P\bar{v}_i = P(I - QP)e_i \stackrel{\eqref{eq:Q}}{=} 0}.$$ 
Then~\eqref{eq:lemma1-linear-1} implies
$$P(S-\Delta L)\bar{v}_i = 0, \qquad i  \in \{1, \ldots, \nu\},$$ 
or, equivalently,
$$
P(S-\Delta L) (I - QP)  =0. 
$$
Finally, by assumption (A$_{Q}$), the matrix $Q$ is a right inverse of $P$ defined as~\eqref{eq:Q}. Then $Q$ is a minimum-norm least-squares solution of the equation ${C\Pi (I - QP) = 0}$~\cite[p.117]{ben2003generalized}. We conclude that all admissible parameters ${P}$, ${\Delta}$ and ${Q}$ yield a solution of the optimization problem~\eqref{eq:optimization-problem-proof} and, thus, the claim is proved.
\end{proof}

\noindent
Theorem~\ref{thm:1} establishes that the  family of models~\eqref{eq:system-rom}-\eqref{eq:family-linear} achieves least squares moment matching (at ${\spectrum{S}}$). The parameterization~\eqref{eq:family-linear} also admits a natural system-theoretic interpretation in terms of a two-step model reduction process, as discussed in the next section.

\subsection{A system-theoretic interpretation} \label{ssec:system-theoretic-interpretation}

\begin{figure}[t!]
\centering
\input{least-squares-moment-matching-linear-1.tex}
\centering
\caption{{Least squares model reduction by moment matching as a two-step model reduction process via surrogate models.}}
\label{fig:least-squares-moment-matching-linear-1}
\end{figure}
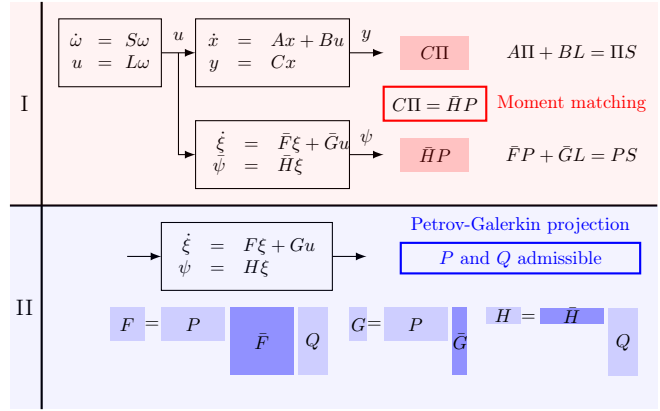

The family of models~\eqref{eq:system-rom}-\eqref{eq:family-linear} can be interpreted in system-theoretic terms as a two-step model reduction process which comprises the following basic steps. 
\begin{itemize}
\item[(I)] The first step consists in constructing a model which incorporates \emph{all} interpolation constraints~\eqref{eq:matching-condition-linear}. This leads to a model of order $\nu$, with $\nu > 2r$. The model is obtained using the parameterization~\eqref{eq:family-linear-astolfi} and is described by the matrices
\beq \label{eq:system-rom-surrogate}
\bar F= S-\Delta L ,
\quad \bar  G=\Delta , 
\quad  \bar  H = C \Pi , 
\eeq
where ${S\in \R^{\nu \times \nu}}$ is a non-derogatory matrix with characteristic polynomial~\eqref{eq:characteristic-polynomial} such that ${\spectrum{A}  \cap \spectrum{S} = \emptyset}$, ${L\in\R^{1\times \nu}}$ is such that the pair $(S,L)$ is observable, ${\Pi \in \R^{n \times \nu}}$ is the (unique) solution of the Sylvester equation~\eqref{eq:Sylvester-equation-astolfi}, and ${\Delta \in \R^{\nu\times 1}}$ is an admissible free parameter, respectively. 
\item[(II)] The second step consists in constructing a reduced order model of~\eqref{eq:system-rom-surrogate} through a biorthogonal Petrov-Galerkin projection, defined as
\beq \label{eq:system-rom-surrogate-reduced}
F=P \bar F Q,
\quad G=P \bar G, 
\quad  H = \bar H Q ,
\eeq
in which ${P \in \R^{r \times \nu}}$ and ${Q \in\R^{\nu \times r}}$  are admissible free parameters.
\end{itemize}
Note that~\eqref{eq:system-rom-surrogate} and~\eqref{eq:system-rom-surrogate-reduced} together yield~\eqref{eq:family-linear}, which shows that the family of models~\eqref{eq:system-rom}-\eqref{eq:family-linear} can be indeed described  as a two-step model reduction process. 
According to this interpretation, the family of models~\eqref{eq:system-rom}-\eqref{eq:family-linear} is obtained by approximation of an auxiliary, possibly large, model~\eqref{eq:system-rom-surrogate}  which takes into account all interpolation constraints~\eqref{eq:matching-condition-linear}. For this reason,~\eqref{eq:system-rom-surrogate} is referred to as a \emph{surrogate model}. Fig.~\ref{fig:least-squares-moment-matching-linear-1} provides a diagrammatic illustration of  of least squares moment matching via surrogate models.

\begin{figure}[t!]
\centering
\input{least-squares-moment-matching-linear-2.tex}
\centering
\caption{{Least squares model reduction by moment matching as a two-step model reduction process via surrogate signal generators.}}
\label{fig:least-squares-moment-matching-linear-2}
\end{figure}
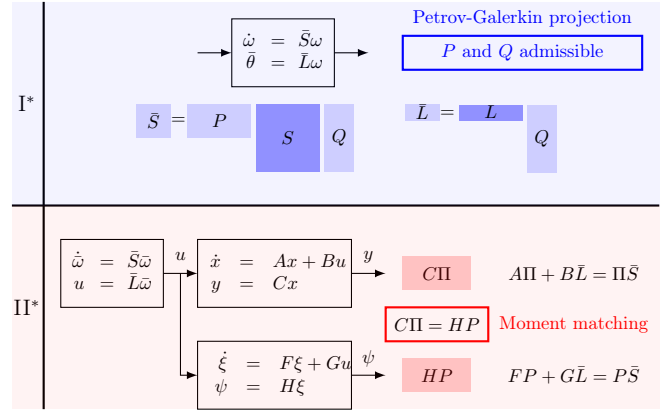

An alternative interpretation of the family of models~\eqref{eq:system-rom}-\eqref{eq:family-linear}  can be given in terms of a ``dual'' two-step model reduction process, which makes use of an auxiliary signal generator.
\begin{itemize}
\item[(I$^*$)] The first step consists in constructing reduced order model of the signal generator~\eqref{eq:system-signal-generator} through a biorthogonal Petrov-Galerkin projection. This step overcomes the issue of taking into account all interpolation constraints~\eqref{eq:matching-condition-linear} by approximating the signal generator~\eqref{eq:system-signal-generator}  with  a ``reduced'' signal generator of order $r$,  with $ r < \nu /2 $. The reduced signal generator in question is described by the matrices
\beq \label{eq:system-signal-generator-surrogate}
\bar S= P S Q ,
\quad \bar  L=  L Q ,
\eeq
in which ${S\in \R^{\nu \times \nu}}$ is a non-derogatory matrix with characteristic polynomial~\eqref{eq:characteristic-polynomial} such that ${\spectrum{A} \cap \spectrum{S}= \emptyset}$, 
${L\in\R^{1\times \nu}}$ is such that the pair $(S,L)$ is observable, ${P \in \R^{r \times \nu}}$  and $Q \in \R^{\nu \times r}$ are admissible free parameters, respectively.

\item[(II$^*$)] The second step consists in building a model of the reduced signal generator using the parameterization~\eqref{eq:family-linear-astolfi}, defined as
\beq \label{eq:system-signal-generator-surrogate-reduced}
F = \bar  S-  \bar \Delta \bar L ,
\quad G=\bar  \Delta , 
\quad  \bar  H = C \bar \Pi .
\eeq 
where $\bar \Pi \in \R^{n \times r }$ is the (unique) solution of the reduced Sylvester equation
\beq \label{eq:Sylvester-equation-surrogate}
A\bar\Pi +B \bar L=\bar \Pi \bar S,
\eeq
and ${\Delta \in \R^{r \times 1}}$ is any admissible vector such that $\spectrum{\bar  S-  \bar \Delta \bar L } \cap  \spectrum{\bar S}  = \emptyset$. Since $\rank P = r$, the vector $\bar \Delta$ can be conveniently rewritten  as $\bar \Delta = P \Delta$, which shows that the family~\eqref{eq:system-signal-generator-surrogate-reduced} can be equivalently defined as
\beq \label{eq:system-signal-generator-surrogate-reduced-equivalent}
F = \bar  S-  P \Delta \bar L ,
\quad G=P \Delta, 
\quad  \bar  H = C \bar \Pi .
\eeq 
\end{itemize}
Note that~\eqref{eq:system-signal-generator-surrogate} and~\eqref{eq:system-signal-generator-surrogate-reduced-equivalent} together yield~\eqref{eq:family-linear}, which shows that the family of models~\eqref{eq:system-rom}-\eqref{eq:family-linear} can be indeed described  as a two-step model reduction process. 
According to this interpretation, the family of models~\eqref{eq:system-rom}-\eqref{eq:family-linear} is obtained by approximation of an auxiliary (possibly large) signal generator~\eqref{eq:system-signal-generator-surrogate} which is then used to construct a reduced order model. For this reason,~\eqref{eq:system-signal-generator-surrogate} is referred to as a \emph{surrogate signal generator}. Fig.~\ref{fig:least-squares-moment-matching-linear-2} provides a diagrammatic illustration of least squares moment matching via surrogate signal generators.

\begin{remark} \label{rem:convex}
The interpretation of least squares moment matching in terms of a two-step model reduction process allows one to compute models achieving least squares moment matching by means of standard convex optimization tools (e.g., CVX~\cite{cvx2014}). As a matter of fact, for every given full rank matrix ${P\in\R^{r\times \nu}}$,  with $\ker P$ a $(S,L)$ conditioned invariant,   the optimization problem~\eqref{eq:optimization-problem} boils down to a constrained least squares problem subject to the additional (generic) constraint $ \spectrum{F} \cap \spectrum{S} = \emptyset$.  Thus, fixing the matrix ${P}$ \textit{a priori}, one can compute a model~\eqref{eq:system-rom} achieving least squares moment matching \textit{directly} by solving the relaxed optimization problem 
\beq \label{eq:optimization-problem-relaxed}
\begin{array}{ll}
    \mbox{minimize}      & \norm{(C \Pi - H P)T }_{2*}^2 \\
    \mbox{subject to}    & FP +GL= P S , 
\end{array}
\eeq
with optimization variables ${F\in\R^{r \times r}}$, ${G\in\R^{r \times 1}}$, ${H\in\R^{1\times r}}$, and verifying that the generic constraint $ \spectrum{F} \cap \spectrum{S} = \emptyset$ holds \textit{a posteriori}. 
We illustrate this point by means of a worked-out example in Section~\ref{sec:examples}.
\end{remark}

\subsection{Least squares moment matching with given properties} \label{ssec:properties}

Model reduction methods are often required to preserve or enforce prescribed properties~\cite{antoulas2005approximation}, including stability, passivity, and a given $L_2$-gain. In practice, this translates into properly selecting the projector $P$ or the vector $\Delta$ which parameterize the families~\eqref{eq:family-linear-astolfi} and~\eqref{eq:family-linear-Krylov} (see, e.g.,~\cite[p.358]{antoulas2005approximation} and~\cite{astolfi2010model}). Similarly, the available degrees of freedom in the parameterization~\eqref{eq:family-linear} can be used to assign specific properties to models achieving least squares moment matching, selecting \textit{both} the projector $P$ and the vector $\Delta$ appropriately. Since least squares moment matching can be interpreted as a two-step model reduction process, one can take advantage of available methods in the literature to assign prescribed properties to the surrogate model/signal generator and, subsequently, to the reduced order model. For reasons of space, we do not provide an exhaustive solution to the general problem of preserving prescribed properties in models achieving least squares moment matching. However, we illustrate how to select the parameters $P$ and $\Delta$  in order to preserve the dominant eigenvalues of a system (and, hence, its stability properties).

\subsubsection{Least squares moment matching with prescribed eigenvalues}

The parameterization~\eqref{eq:family-linear} can be used to preserve the first ${r}$ dominant eigenvalues of a system by selecting the parameters $P$ and $\Delta$ as follows. First, the vector $\Delta$ is selected in such a way that the matrix $(S-\Delta L)$ preserves the first ${\nu}$ dominant eigenvalues of the original system. This amounts to assigning the spectrum of the surrogate model~\eqref{eq:system-rom-surrogate}, which is always possible since the pair $(S,L)$ is observable, by assumption.  Second, the matrix $P$ is defined as 
$P = [\, P_1^{\transpose} \ \cdots \ P_{r}^{\transpose} \,]^{\transpose} ,$
with $\{P_1,\ldots, P_r\}$ a real Jordan basis of right eigenvectors of the matrix $S-\Delta L$ corresponding to the first ${r}$ dominant eigenvalues of the matrix ${(S-\Delta L)}$. By construction, this yields 
\beq \label{eq:properties}
P (S-\Delta L) = F P
\eeq
for some matrix ${F \in \R^{ r \times r}}$. This ensures that the parameters $P$ and $\Delta$ are admissible, since $P$ is full rank and such that the subspace $\ker P$ is ${(S-\Delta L)}$-invariant (and, thus, also $(S,L)$ conditioned invariant). Furthermore, \eqref{eq:properties} ensures that the reduced order model preserves the first ${r}$ dominant eigenvalues of the matrix ${(S-\Delta L)}$ and, hence, those of the original system. We illustrate this construction with a worked-out example in Section~\ref{sec:examples}.

\section{Nonlinear systems} \label{sec:nonlinear}

This section provides a nonlinear enhancement of the notion of least squares moment matching. 
In close analogy with Section~\ref{sec:linear}, the notion of least squares moment matching is then used to develop a model reduction theory for nonlinear systems. 
Our analysis builds on the formalism presented in~\cite{astolfi2010model}, which we first briefly recall for completeness, with minor variations.

\subsection{Moment matching}

Consider a continuous-time, single-input, single-output, time-invariant system described by the equations\footnote{Throughout the paper all mappings are assumed to be smooth, \emph{i.e.} infinitely many times differentiable, if not otherwise stated.}
\beq \label{eq:system-nonlinear}
\dot{x} = f(x,u), \quad y=h(x),
\eeq 
with ${x(t) \in\R^{n}}$, ${u(t)\in\R}$, ${y(t)\in\R}$, and ${f:\R^n \times \R \to\R^n}$, ${g:\R^n \to\R^n}$ and ${h:\R^n\to\R}$ such that ${f(0,0)=0}$ and ${h(0)=0}$, and a signal generator described by the equations
\beq \label{eq:signal-generator-nonlinear}
\dot{\omega} = s(\omega), \quad \theta =l(\omega),
\eeq 
in which ${\omega(t) \in \Omega}$ and ${\theta(t)\in\R}$, with ${\Omega \subset \R^{\nu}}$ a sufficiently small\footnote{All statements are local, although global versions can be easily given.} open, connected, invariant neighborhood containing the origin, while the mappings ${s:\Omega \to \Omega}$ and ${h:\Omega\to\R}$ are such that ${s(0)=0}$ and ${l(0)=0}$, respectively.

For the notion of moment to make sense, we make the following standing assumption throughout the paper.

\begin{standing assumption}
The system~\eqref{eq:system-nonlinear} is minimal, \textit{i.e.} (locally) accessible and observable\footnote{See~\cite[Section 3]{nijmeijer1990nonlinear} for the notion of local accessibility and observability. See also~\cite[Section 2.2]{scarciotti2017nonlinear} for further detail on the assumption of minimality in the context of moment matching for nonlinear systems.} at the origin.
\end{standing assumption}

\noindent
The definition of moment of a nonlinear system also requires the following assumptions.

\begin{assumption} \label{ass:signal-generator-nonlinear}
The signal generator~\eqref{eq:signal-generator-nonlinear} is (locally) observable at the origin
and such that the matrix  ${S = \tfrac{\partial s}{\partial \omega}{\scriptstyle(0)}}$ is non-derogatory.
\end{assumption}

\begin{assumption} \label{ass:moment-well-defined-nonlinear}
The partial differential equation
\beq \label{eq:moment-matching-PDE-system}
f({\pi}(\omega),l(\omega)) = \frac{\partial {\pi} }{\partial \omega} (\omega) s(\omega),
\eeq
admits a unique solution $\pi(\cdot)$, locally defined in the neighbourhood ${\Omega}$ of the origin, such that ${\pi(0) =0}$. 
\end{assumption} 

\begin{definition}~\cite{astolfi2010model}
Consider system~\eqref{eq:system-nonlinear} and the signal generator~\eqref{eq:signal-generator-nonlinear}. 
Suppose Assumptions~\ref{ass:signal-generator-nonlinear} and~\ref{ass:moment-well-defined-nonlinear} hold. The \emph{moment of system~\eqref{eq:system-nonlinear} at $(s,l)$} is defined as the mapping  ${\mu(\cdot) = h (\pi(\cdot))},$  where ${\pi(\cdot)}$ is the unique solution of the partial differential equation~\eqref{eq:moment-matching-PDE-system}.
\end{definition}

\begin{remark}
Assumptions~\ref{ass:signal-generator-nonlinear} and~\ref{ass:moment-well-defined-nonlinear} ensure that the moment of system~\eqref{eq:system-nonlinear} at $(s,l)$ is (locally) well-defined, thus providing a direct nonlinear counterpart of the assumptions made on the signal generator~\eqref{eq:system-signal-generator} in the linear case. We emphasize that Assumption~\ref{ass:signal-generator-nonlinear} is a strong
observability assumption, which directly implies the observability assumption originally used in~\cite{astolfi2010model} and which we make to simplify the exposition. We also emphasize that the partial differential equation~\eqref{eq:moment-matching-PDE-system} admits a unique \textit{formal} (power series) solution if and only if the following non-resonance condition holds~\cite[Lemma 4.13]{huang2004nonlinear}
\beq \label{eq:non-resonance-condition}
\spectrum{A} \cap \spectrumk{S} = \emptyset, \quad k = 1, 2, \dots ,
\eeq
in which ${A = \tfrac{\partial f}{\partial x}(0,0)}$ and ${S = \tfrac{\partial s}{\partial \omega}(0)}$, and 
\beq \label{eq:spectrum-extended}
\spectrumk{S} = \left\{\,\lambda \in \Co \,:\, \lambda = \sum_{i = 1}^\nu \lambda_i k_i , \ k = \sum_{i = 1}^\nu k_i   \, \right\},
\eeq
where ${\lambda_i \in \spectrum{S}}$ and ${k_i \in \{ 0,1,\dots, k\}}$. 
Assumption~\ref{ass:moment-well-defined-nonlinear} therefore holds if the equilibrium ${x=0}$ of
the system ${\dot{x} = f(x,0)}$ is hyperbolic (or, in particular, locally exponentially stable) and if the signal generator~\eqref{eq:signal-generator-nonlinear} is Poisson stable\footnote{See \cite[Chapter 8]{isidori1995nonlinear} for the definition of Poisson stability.} (or, in particular, periodic). For real analytic systems~\eqref{eq:system-nonlinear} and~\eqref{eq:signal-generator-nonlinear}, the non-resonance condition~\eqref{eq:non-resonance-condition} also implies that the unique formal power series solution converges to a real analytic mapping~\cite[Theorem 1]{kazantzis2000singular}, provided that zero is not in the convex hull of $\spectrum{S}$ (see also~\cite{krener2002nonlinear,krener2004erratum} for more detail). 
\end{remark}

\noindent
In analogy with the linear case, it is possible to establish a one-to-one correspondence between the moment of system~\eqref{eq:system-nonlinear} at $(s,l)$ and  the steady-state output response of the interconnected system~\eqref{eq:system-nonlinear}-\eqref{eq:signal-generator-nonlinear}, with ${u =  \theta}$, provided that the following excitability assumption holds.

\begin{assumption} \label{ass:excitability-nonlinear}
The vector ${\omega(0) \in\R^{\nu}}$ is such that the pair $(s,\omega(0))$ is exciting\footnote{See \cite{padoan2016geometric} for the definition exciting pair.}.
\end{assumption}

\begin{theorem}~\cite{astolfi2010model} \label{thm-astolfi-1-nonlinear} 
Consider system~\eqref{eq:system-nonlinear} and the signal generator~\eqref{eq:signal-generator-nonlinear}. 
Suppose Assumptions~\ref{ass:signal-generator-nonlinear},~\ref{ass:moment-well-defined-nonlinear} and~\ref{ass:excitability-nonlinear} hold. Assume that the zero equilibrium of the system ${\dot{x} = f(x,0)}$ is locally exponentially stable and that system~\eqref{eq:signal-generator-nonlinear} is Poisson stable.  Then the moment of system~\eqref{eq:system-nonlinear} at $(s,l)$ is (locally) well-defined and uniquely determined by the steady-state output response\footnote{See \cite[Chapter 8]{isidori1995nonlinear} for the notion of steady-state response.} of the interconnected   
system~\eqref{eq:system-nonlinear}-\eqref{eq:signal-generator-nonlinear}, with ${u =  \theta}$.
\end{theorem}

\noindent
Theorem~\ref{thm-astolfi-1-nonlinear} motivates the definition of a nonlinear counterpart of the notion of model at $(S,L)$ for a system described by the equations
\beq \label{eq:system-model-nonlinear}
\dot{\xi} = \phi(\xi,v) , \quad  \psi = \kappa(\xi),
\eeq
with ${\xi(t)\in\R^{r}}$, ${v(t)\in\R}$, ${\psi(t)\in\R}$, and ${\phi:\R^r \times \R \to\R^n}$, ${\kappa:\R^r \to\R}$  such that ${\phi(0,0)=0}$ and ${\kappa(0)=0}$, respectively.

\begin{definition}~\cite{astolfi2010model}
The system~\eqref{eq:system-model-nonlinear} is a \emph{model of system~\eqref{eq:system-nonlinear} at $(s,l)$}  if its moment at $(s,l)$ is (locally) well-defined and coincides with that of system~\eqref{eq:system-nonlinear}, \emph{i.e.} if the partial differential equation
\beq \label{eq:moment-matching-PDE-model}
\phi(p(\omega),l(\omega)) = \frac{\partial p }{\partial \omega}(\omega) s(\omega)
\eeq
possesses a unique solution ${p(\cdot)}$, locally defined in the neighbourhood $\Omega$ of the origin, such that ${p(0)=0}$ and
\beq \label{eq:moment-matching-nonlinear}
h(\pi(\omega)) = \kappa(p(\omega)),
\eeq
where ${\pi(\cdot)}$ is the (unique) solution of the partial differential equation~\eqref{eq:moment-matching-PDE-system}. In this case, system~\eqref{eq:system-model-nonlinear} is said to \emph{match the moment of system~\eqref{eq:system-nonlinear} (or to achieve moment matching) at $(s,l)$}. Furthermore, system~\eqref{eq:system-model-nonlinear} is a \emph{reduced order model of system~\eqref{eq:system-nonlinear} at $(s,l)$} if ${r<n}$.
\end{definition}

\noindent
A family of reduced order models achieving moment matching for system~\eqref{eq:system-nonlinear} has been determined in~\cite{astolfi2010model} by selecting $r=\nu$ and ${p(\omega)=\omega}$ in~\eqref{eq:moment-matching-PDE-model} and~\eqref{eq:moment-matching-nonlinear}, respectively. As a matter of fact, this yields a family of reduced order models of system~\eqref{eq:system-nonlinear} at $(s,l)$ described by the equations
\beq \label{eq:family-nonlinear-astolfi}
\dot{\xi} = s(\xi) -\delta(\xi) l(\xi)  + \delta(\xi) u, \quad  \psi = h(\pi(\xi)),
\eeq
with ${\delta: \R^{r} \to \R^{r}}$ any mapping such that the equation
\beq
s(p(\omega)) -\delta(p(\omega)) l(p(\omega))  + \delta(p(\omega)) l(\omega) = \frac{\partial p }{\partial \omega} s(\omega)
\eeq
has the unique solution ${p(\omega)=\omega}$. In analogy with the linear case, the mapping ${\delta(\cdot)}$ is a \emph{parameter} of the family of reduced order models~\eqref{eq:family-nonlinear-astolfi} which can be used to assign prescribed properties to the reduced order model, including stability, passivity, and a given $L_2$-gain~\cite{astolfi2010model}.

\subsection{Least squares moment matching}

This section introduces the notion of least squares moment matching for nonlinear systems, motivated by the equivalence established in Theorem~\ref{thm:moments-optimization-linear}. To this end, we introduce the following assumption.

\begin{assumption} \label{ass:T-nonlinear}
The mapping ${\tau: \Omega \to \Omega}$ is a diffeomorphism and such that ${\tau(0)=0}$. 
\end{assumption}

\begin{definition}  \label{def:least-squares-moment-matching-nonlinear}
Consider system~\eqref{eq:system-nonlinear} and the signal generator~\eqref{eq:signal-generator-nonlinear}. 
Suppose Assumptions~\ref{ass:signal-generator-nonlinear},~\ref{ass:moment-well-defined-nonlinear} and \ref{ass:T-nonlinear}   hold.
The model~\eqref{eq:system-model-nonlinear} \emph{achieves least squares moment matching at $(s,l)$} if the triple $(\phi,\kappa,p)$ is a (formal) solution of the constrained the optimization problem 
\beq \label{eq:optimization-problem-nonlinear}
\begin{array}{ll}
    \mbox{minimize}      
    &\displaystyle\sup_{ \omega \in \Omega}\norm{\left(\frac{\partial \mu }{\partial \omega}(\tau(\omega))-\frac{\partial \hat{\mu} }{\partial \omega} (\tau(\omega))\right)\frac{\partial \tau }{\partial \omega}(\omega)}_{2*} 
    \\[0.5cm]
    \mbox{subject to}    
    & \phi(p(\omega),l(\omega)) = \frac{\partial p }{\partial \omega}(\omega) s(\omega), \ \omega \in \Omega, 
    \\[0.3cm]
    & \spectrum{F} \cap \spectrumk{S} = \emptyset, \  k = 1, 2, \dots , 
\end{array}
\eeq
in which ${\mu(\cdot) = h(\pi(\cdot))}$, ${\hat\mu(\cdot) = \kappa(p(\cdot))}$, ${F = \tfrac{\partial \phi}{\partial \xi}{\scriptstyle(0,0)}}$ and ${S = \tfrac{\partial s}{\partial \omega}{\scriptstyle(0)}}$,  where ${\phi(\cdot,\cdot)}$, ${\kappa(\cdot)}$ and ${p(\cdot)}$ are the optimization variables, while system~\eqref{eq:system-nonlinear} and the signal generator~\eqref{eq:signal-generator-nonlinear} (and, thus, the solution ${\pi(\cdot)}$ of the partial differential equation~\eqref{eq:moment-matching-PDE-system})  are problem data. 
\end{definition}

We are now in a position to pose the \textit{model reduction problem by least squares moment matching} for nonlinear systems.

\begin{problem}
Consider system~\eqref{eq:system-nonlinear} and the signal generator~\eqref{eq:signal-generator-nonlinear}. 
Suppose Assumptions~\ref{ass:signal-generator-nonlinear},~\ref{ass:moment-well-defined-nonlinear} and \ref{ass:T-nonlinear}  hold. Let ${r\in \Zplus}$, with ${2r < \nu}$.  Find, if possible, a system~\eqref{eq:system-model-nonlinear} of order $r$ which achieves least squares moment matching at $(s,l)$.
\end{problem}

Before proceeding to the solution of the model reduction problem by least squares moment matching for nonlinear systems, we show that models achieving least squares moment matching at $(s,l)$ directly provide an \textit{a priori} error bound in terms of the worst case r.m.s. gain of an error system with respect to the family of signals produced by the signal generator~\eqref{eq:signal-generator-nonlinear}, in close analogy with the linear case.

\subsubsection{An \textit{a priori}  error bound}

We now characterize least squares moment matching in terms of the steady-state behavior of the error system 
\beq \label{eq:system-error-nonlinear}
\dot{x} = f(x,u), \quad
\dot{\xi} = \phi(\xi,u) , \quad  
 e=h(x)-\kappa(\xi),  
\eeq 
in which  ${x(t)\in\R^n}$, ${\xi(t)\in\R^r}$, ${u(t)\in\R}$, and ${e(t)\in\R}$.
Generalizing Theorem~\ref{thm:rms-linear}, we show that if the steady-state output response $e_{ss}$ of the interconnected system~\eqref{eq:signal-generator-nonlinear}-\eqref{eq:system-error-nonlinear}, with ${u = \theta}$, is (locally) well-defined and if all solutions of the signal generator~\eqref{eq:system-signal-generator} are periodic, then achieving least squares moment matching corresponds to minimizing an upper bound of the 
\textit{worst case r.m.s. gain} of the error system~\eqref{eq:system-error-nonlinear}  with respect to the family of  signals produced by the signal generator~\eqref{eq:system-signal-generator}, defined as\footnote{See~\cite{isidori1992disturbance} for closely the related notion of $H_{\infty}$ gain of a nonlinear system.}
\beq \label{eq:gain-rms-nonlinear}
\gamma_{rms} =  \sup_{\omega(\cdot) \in \mathcal{W}} \frac{\norm{e_{ss}}_{rms}}{\norm{\omega}_{rms}} ,
\eeq 
where $\mathcal{\mathcal{W}}$ is the family of signals produced by~\eqref{eq:signal-generator-nonlinear} with non-zero r.m.s. value.

\begin{theorem} \label{thm:rms-nonlinear}
Consider system~\eqref{eq:system-nonlinear}, the signal generator~\eqref{eq:signal-generator-nonlinear} and the model~\eqref{eq:system-model-nonlinear}. Suppose Assumptions~\ref{ass:moment-well-defined-nonlinear},~\ref{ass:excitability-nonlinear} and~\ref{ass:T-nonlinear} hold. Assume that the zero equilibrium of the system ${\dot{x} = f(x,0)}, \, {\dot{\xi} = \phi(\xi,0)}$ is locally exponentially stable and that all solutions of system~\eqref{eq:signal-generator-nonlinear} are periodic. Then the following statements hold.
\begin{itemize}
\item[(i)] The steady-state output response of the interconnected system~\eqref{eq:signal-generator-nonlinear}-\eqref{eq:system-error-nonlinear}, with ${u = \theta}$, is (locally) well-defined and uniquely determined by the moment of the error system~\eqref{eq:system-error-nonlinear} at $(s,l)$.
\item[(ii)]  The worst case r.m.s. gain of the error system~\eqref{eq:system-error-nonlinear}  with respect to the family of signals produced by the signal generator~\eqref{eq:system-signal-generator} is well-defined and such that
\beq \label{eq:gain-bound-nonlinear}
\gamma_{rms}  \le \displaystyle\sup_{ \omega \in  \Omega}\norm{\frac{\partial \mu }{\partial \omega}(\omega)-\frac{\partial \hat{\mu} }{\partial \omega} (\omega)}_{2*}  ,
\eeq
where ${\mu(\cdot)}$ and ${\hat{\mu}(\cdot)}$ are the moments of systems~\eqref{eq:system-nonlinear} and~\eqref{eq:system-model-nonlinear} at $(s,l)$, respectively.
\item[(iii)]   The error bound~\eqref{eq:gain-bound-nonlinear} is minimized if the model~\eqref{eq:system-model-nonlinear} achieves least squares moment matching at~$(s,l)$ and if the mapping $\tau(\cdot)$ is an isometry\footnote{See~\cite[p.332]{lee2013smooth} for the definition.} (on $\Omega$). 
\end{itemize}
\end{theorem}

\begin{proof}
(i). The dynamics of the interconnected system~\eqref{eq:signal-generator-nonlinear}-\eqref{eq:system-error-nonlinear}, with ${u = \theta}$, is governed by the equations
\beq \label{eq:system-error-interconnected-nonlinear} 
\dot{\omega} = s(\omega), \
\dot{x} = f(x,l(\omega)),  \
\dot{\xi} = \phi(\xi,l(\omega)),   \ 
e=h(x)-\kappa(\xi). 
\eeq 
By assumption, the equilibrium at the origin of the system ${\dot{x} = f(x,0), \ \dot{\xi} = \phi(\xi,0)}$ is locally exponentially stable and the system~\eqref{eq:signal-generator-nonlinear} is periodic. By the center manifold theorem~\cite[p.4]{carr1981application}, the interconnected system~\eqref{eq:system-error-interconnected-nonlinear} has therefore a well-defined center manifold
\beq \nn
\mathcal{M}_c = \left\{ \, (x,\xi,\omega)  \,:\,  x = \pi(\omega), \, \xi = p(\omega) \, \right\} ,
\eeq
where ${\pi(\cdot)}$ and ${p(\cdot)}$ are (locally) the unique solutions of the partial differential equations~\eqref{eq:moment-matching-PDE-system} and~\eqref{eq:moment-matching-PDE-model}, respectively. Furthermore, the manifold $\mathcal{M}_c$ is (locally) exponentially attractive, since
\beq\nn
\renewcommand{\arraystretch}{.9}
\setlength{\arraycolsep}{1pt}
\dot{\overbrace{\bma \begin{array}{c} x- \pi(\omega) \\ \xi - p(\omega) \end{array} \ema}} 
\stackrel{\eqref{eq:system-error-interconnected-nonlinear}}{=} 
	\bma 
		\begin{array}{c} 
			f(x,l(\omega)) - \frac{\partial\pi}{\partial \omega} s(\omega) \\ 
			\phi(\xi,l(\omega))  - \frac{\partial p}{\partial \omega} s(\omega)
		\end{array} 
	\ema  
\stackrel{\eqref{eq:moment-matching-PDE-system}, \eqref{eq:moment-matching-PDE-model}}{=} 
	\bma 
		\begin{array}{c} 
			f(x,l(\omega)) - f(\pi(\omega),l(\omega)) \\ 
			\phi(\xi,l(\omega))  - \phi(\pi(\omega),l(\omega)) 
		\end{array} 
	\ema .
\eeq
The output of the interconnected system~\eqref{eq:system-error-interconnected-nonlinear}
can be thus written as ${e =  e_{ss} + e_{d}}$ in which
\beq \label{eq:error-steady-state-nonlinear}
e_{ss} = \mu (\omega) - \hat \mu (\omega) = h(\pi(\omega)) - \kappa(p(\omega)) ,
\eeq
is the (locally) well-defined steady-state output response of the interconnected system~\eqref{eq:system-error-interconnected-nonlinear} and $e_{d}$ is an exponentially decaying signal. Furthermore, under the stated assumptions, Theorem~\ref{thm-astolfi-1-nonlinear} implies that the steady-state output response of the interconnected system~\eqref{eq:system-error-interconnected-nonlinear} is uniquely determined by the moment of the error system~\eqref{eq:system-error-nonlinear} at $(s,l)$.

(ii). The worst case r.m.s. gain of the error system~\eqref{eq:system-error-nonlinear}  with respect to the family of  signals produced by the signal generator~\eqref{eq:system-signal-generator} is well-defined: by assumption, all solutions of the signal generator~\eqref{eq:signal-generator-nonlinear} are periodic and, thus, the steady-state output response of the interconnected system~\eqref{eq:system-error-interconnected-nonlinear} is also periodic (with the same period) by~\eqref{eq:error-steady-state-nonlinear}. To obtain the error bound~\eqref{eq:gain-bound-nonlinear}, 
note that
\begin{align}
\!\norm{e_{ss}}_{2} \!\!\!\!\!
&~~~\stackrel{\eqref{eq:error-steady-state-nonlinear}}{=}~
\norm{\mu(\omega) - \hat{\mu}(\omega)}_{2}  \nn \\
&~~~~~{\le}~
\norm{\left(\int_{0}^1 
\frac{\partial \mu}{\partial \omega} (\alpha \omega) - 
\frac{\partial \hat \mu}{\partial \omega} (\alpha \omega)
d\alpha \right) \omega }_{2}\nn\\
&~~~~~{\le}~
\norm{\int_{0}^1 
\frac{\partial \mu}{\partial \omega} (\alpha \omega) - 
\frac{\partial \hat \mu}{\partial \omega} (\alpha \omega)
d\alpha }_{2*} \! \! \norm{ \omega }_{2}\nn\\
&~~~~~{\le}~
\left(\int_{0}^1 \norm{ 
\frac{\partial \mu}{\partial \omega} (\alpha \omega) - 
\frac{\partial \hat \mu}{\partial \omega} (\alpha \omega)
}_{2*} \!\!\!\! d\alpha \right) \norm{ \omega }_{2}\nn\\
&~~~~~{\le}~
\left(\int_{0}^1 \sup_{\omega \in \Omega} \norm{  
\frac{\partial \mu}{\partial \omega} (\omega) - 
\frac{\partial \hat \mu}{\partial \omega} (\omega)
}_{2*}  \! d\alpha \right) \norm{ \omega }_{2}\nn\\
&~~~~\,=~
\left( \sup_{\omega \in \Omega} 
\norm{\frac{\partial \mu }{\partial \omega}(\omega)-\frac{\partial \hat{\mu} }{\partial \omega} (\omega)}_{2*} 
\right) \norm{ \omega }_{2} , \label{eq:coordinates-change-nonlinear}
\end{align} 
where the first inequality follows from the mean value theorem~\cite[p.40]{nijmeijer1990nonlinear} 
and the second inequality follows from the Cauchy-Schwartz inequality.
Then
\beq \nn
\gamma_{rms} 
\!\stackrel{\eqref{eq:gain-rms-nonlinear}}{=} \! 
\sup_{\omega(\cdot) \in \mathcal{W}} \frac{\norm{e_{ss}}_{rms}}{\norm{\omega}_{rms}} 
\!\stackrel{\eqref{eq:coordinates-change-nonlinear}}{\le} 
\sup_{\omega \in \Omega} 
\norm{\frac{\partial \mu }{\partial \omega}(\omega)-\frac{\partial \hat{\mu} }{\partial \omega} (\omega)}_{2*},
\eeq
which proves the error bound~\eqref{eq:gain-bound-nonlinear}.

(iii) If the mapping $\tau(\cdot)$ is an isometry on $\Omega$, then its differential (expressed in local coordinates by the matrix) $\frac{\partial \tau}{\partial \omega} {\scriptstyle(\omega)}$ is orthogonal for every 
${\omega \in \Omega}$~\cite[p.224]{robbin2021introduction}.
Then
\beq \nn
\norm{
	\frac{\partial \mu }{\partial \omega}{(\tau(\omega))} 
   -\frac{\partial \hat{\mu} }{\partial \omega}{(\tau(\omega))}
	}_{2*} 
\!=\!
\norm{\left(\frac{\partial \mu }{\partial \omega}{(\tau(\omega))} )-\frac{\partial \hat{\mu} }{\partial \omega} {(\tau(\omega))} \right)\frac{\partial \tau }{\partial \omega}{(\omega)}}_{2*} \!,
\eeq
since the norm $\norm{\,\cdot\,}_{2*}$ is unitarily invariant. The claim is thus a direct consequence of item (ii), Definition~\ref{def:least-squares-moment-matching-nonlinear} and the fact that $\tau(\cdot)$ is, by assumption, a diffeomorphism (on $\Omega$).
\end{proof}

\begin{remark}
For linear systems, least squares moment matching does not require any stability assumption on the model, only that the non-resonance condition ${\spectrum{F} \cap \spectrum{S} = \emptyset}$ is satisfied. Similarly, for nonlinear systems, the notion of least squares moment matching is well-defined for any  model~\eqref{eq:system-model-nonlinear} and any signal generator~\eqref{eq:signal-generator-nonlinear} such that the non-resonance condition ${\spectrum{F} \cap \spectrumk{S} = \emptyset}$, with ${F = \tfrac{\partial \phi}{\partial \xi}{\scriptstyle(0,0)}}$ and ${S = \tfrac{\partial s}{\partial \omega}{\scriptstyle(0)}}$, holds for every integer ${k>0}$. However, when considering unstable systems the relation between least squares moment matching and the steady-state response of the error system~\eqref{eq:system-error-nonlinear} established by Theorem~\ref{thm:rms-nonlinear} is lost.
\end{remark}

\subsection{Model reduction by least squares moment matching}

\noindent
The construction of models achieving moment matching requires the determination of mappings $\phi(\cdot,\cdot)$, $\kappa(\cdot)$ and $p(\cdot)$ which (formally) solve the optimization problem~\eqref{eq:optimization-problem-nonlinear} in a neighbourhood of the origin.
Finding analytical expressions for such mappings is generally hard.
However, in analogy with the linear case, 
we circumvent this issue by regarding the optimization variable $p(\cdot)$ as a given \emph{parameter} and introducing a nonlinear counterpart of the family of models~\eqref{eq:system-rom}-\eqref{eq:family-linear}.

Consider system~\eqref{eq:system-nonlinear} and the signal generator~\eqref{eq:signal-generator-nonlinear}. Suppose Assumptions~\ref{ass:signal-generator-nonlinear} and~\ref{ass:moment-well-defined-nonlinear} hold. The family of models in question is described by the equations
\beq \label{eq:family-nonlinear}
\dot{\xi} = \varphi(\xi) + \gamma (\xi)u , \quad  \psi = \kappa(\xi),
\eeq
with
\bseq  \label{eq:family-nonlinear-varphi}
\begin{align}
\varphi(\xi)	&= 	\left.\frac{\partial p}{\partial \omega}(\omega) (s(\omega) - \delta(\omega) l(\omega)) \right|_{\omega = q(\xi)}    \\
\gamma(\xi)	&= 	\left.\frac{\partial p}{\partial \omega}(\omega) \delta (\omega) \right|_{\omega = q(\xi)}, \\
\kappa(\xi)	&= 	\left. h(\pi(\omega)) \right|_{\omega = q(\xi)}, 
\end{align}
\eseq
where ${\pi(\cdot)}$ is the (unique) solution of the partial differential equation~\eqref{eq:moment-matching-PDE-system}, locally defined in the neighbourhood $\Omega$ of the origin and such that $\pi(0)=0$,  while 
${p: \Omega \to \R^{r}}$, ${\delta: \Omega \to \R^{r}}$, and ${q: \R^r \to \Omega}$, with ${p(0)=0}$ and ${q(0)=0}$, are such that 
\begin{itemize}
\item[(A$_p$)]  the mapping $p(\cdot)$ is a surjective submersion\footnote{If $\mathcal{M}$ and $\mathcal{N}$ are smooth manifolds, a smooth map $p: \mathcal{M} \to \mathcal{N}$ is called a \textit{(smooth) submersion} if its differential is surjective at each point~\cite[p.78]{lee2013smooth}.
}
and such that the distribution\footnote{The differential of a smooth mapping $p(\cdot)$ is denoted by $dp(\cdot)$~\cite[p.7]{isidori1995nonlinear}.} $\ker dp$ is $(s,l)$-invariant\footnote{See~\cite{krener1986conditioned} for the definition of $(f,h)$-invariant distribution.},
\item[(A$_q$)] the mapping $q(\cdot)$ is a right inverse of $p(\cdot)$, \textit{i.e.} 
${p(q(\xi)) =\xi}$ for all ${\xi \in \R^r}$, 
defined as
\beq \label{eq:condition-p-nonlinear-0}
q(\xi) = p_{\transpose} \circ (p \circ p_{\transpose} )^{-1} (\xi),
\eeq
in which ${p_{\transpose}(\cdot)}$ is any mapping such that ${p_{\transpose}(0)=0}$ and 
\beq \label{eq:condition-p-nonlinear-2}
\frac{\partial p_{\transpose}}{\partial \omega}(\omega) = 
\frac{\partial \tau}{\partial \omega} (\omega)
\left(\frac{\partial p}{\partial \omega} (\omega) 
\frac{\partial \tau}{\partial \omega} (\omega)
\right)^{\transpose},
\eeq
with ${\tau(\cdot)}$ any diffeomorphism (on $\Omega$) such that ${\tau(0)=0}$,
\item[(A$_\delta$)]  the mapping $\delta(\cdot)$ is such that the distribution $\ker dp$ is invariant\footnote{A distribution $\mathcal{D}$ is invariant under the dynamics of the system ${\dot{x} = f(x)+g(x)u}$, with ${x(t) \in\R^{n}}$ and ${u(t)\in\R}$, if ${[f,\mathcal{D}] \subset \mathcal{D}}$ and  ${[g,\mathcal{D}] \subset \mathcal{D}}$~\cite[Definition 3.44]{nijmeijer1990nonlinear}.} under the dynamics of the system
\beq \label{eq:system-generator-surrogate-nonlinear}
\dot{\omega} = s(\omega)-\delta(\omega) l(\omega) + \delta(\omega) v, 
\eeq
in which ${\omega(t)\in\R^{\nu}}$ and ${v(t)\in\R}$, and such that the non-resonance condition
\beq \label{eq:spectrum-constraint-nonlinear}
{\spectrum{F} \cap \spectrumk{S}  = \emptyset}, \qquad k = 1, 2, \dots ,  
\eeq
holds, with ${F = \tfrac{\partial \varphi}{\partial \xi}{\scriptstyle(0)}}$ and ${S = \tfrac{\partial s}{\partial \omega}{\scriptstyle(0)}}$,
\end{itemize}
in which case $p(\cdot)$, $\delta(\cdot)$ and $q(\cdot)$ are said to be \emph{admissible}.

The family of models~\eqref{eq:family-nonlinear}-\eqref{eq:family-nonlinear-varphi} has the same structure of the family~\eqref{eq:system-rom}-\eqref{eq:family-linear} and, similar to the linear case, provides a parameterization of models achieving least squares moment matching, as detailed by the following statement.    

\begin{theorem} \label{thm:1-nonlinear}
Consider system~\eqref{eq:system-nonlinear} and the family of models~\eqref{eq:family-nonlinear}-\eqref{eq:family-nonlinear-varphi}.
Let $p(\cdot)$, $\delta(\cdot)$ and $q(\cdot)$ be {admissible} for the parameterization~\eqref{eq:family-nonlinear-varphi}. Then the model~\eqref{eq:family-nonlinear}-\eqref{eq:family-nonlinear-varphi} achieves least squares moment matching at ${(s,l)}$.
\end{theorem}

\begin{proof}
The proof closely parallels that of Theorem~\ref{thm:1}, leveraging on tools from nonlinear geometric control theory. The notation is mostly borrowed from~\cite{isidori1995nonlinear} and~\cite{nijmeijer1990nonlinear}. To streamline the exposition, the arguments of the functions are often omitted.

By Definition~\ref{def:least-squares-moment-matching-nonlinear}, we need to show that the triple $(\varphi(\xi)+\gamma(\xi)u, \kappa(\xi), p(\omega))$ is a (formal) solution of the optimization problem~\eqref{eq:optimization-problem-nonlinear} for some (local) diffeomorphism $\tau(\cdot)$ such that ${\tau(0)=0}$.
Equivalently, after simple manipulations, direct substitution of~\eqref{eq:family-nonlinear} and~\eqref{eq:family-nonlinear-varphi} into~\eqref{eq:optimization-problem-nonlinear} shows that it suffices to establish that all admissible parameters yield a solution of the optimization problem
\beq \label{eq:optimization-problem-proof-nonlinear}
\begin{array}{ll}
    \mbox{minimize}      
    & \displaystyle\sup_{\omega \in \Omega}
    		\norm{
    			\frac{\partial \mu}{\partial \omega}
    				\left(
    					I-\frac{\partial q}{\partial \xi}\frac{\partial p}{\partial \omega}
    				\right)
    			\frac{\partial \tau}{\partial \omega}}_{2*},
\end{array}
\eeq
with optimization variables $p(\cdot)$, $q(\cdot)$ and $\delta(\cdot)$ subject to the constraints~\eqref{eq:spectrum-constraint-nonlinear} and 
\beq \label{eq:PDE-nonlinear-admissible}
\frac{\partial p}{\partial \omega}(\omega) s(\omega)
=
\left.\frac{\partial p}{\partial w}(w) (s(w) - \delta(w) l(w)) \right|_{w = q \circ p(\omega)}  +
\left.\frac{\partial p}{\partial w}(w) \delta (w) \right|_{w = q \circ p(\omega)} l(\omega)  
\eeq
for some diffeomorphism $\tau(\cdot)$ (on $\Omega$) such that ${\tau(0)=0}$.

For, consider preliminarily the distribution
${\mathcal{D} = \ker dp}$
and the vector fields 
\beq \label{eq:vector-fields-nonlinear}
\bar{\varphi}(\omega) = s(\omega)-\delta(\omega) l(\omega), \quad \bar{\gamma}(\omega) = \delta(\omega),   \quad \omega \in \Omega .
\eeq
Let $p(\cdot)$, $\delta(\cdot)$ and $q(\cdot)$ be {admissible} for the parameterization~\eqref{eq:family-nonlinear-varphi} and let $\tau(\cdot)$ be a diffeomorphism (on $\Omega$) such that ${\tau(0)=0}$. Note that (A$_\delta$) trivially implies the condition~\eqref{eq:spectrum-constraint-nonlinear}. Furthermore, (A$_p$) and (A$_\delta$) imply that the distribution $\mathcal{D}$ is invariant under the dynamics of system~\eqref{eq:system-generator-surrogate-nonlinear}, \textit{i.e.}  
\begin{align}
&[\bar{\varphi}, \mathcal{D}] \subset  \mathcal{D}, \quad \forall \ V \in \mathcal{D}, \label{eq:invariance-1}\\
&[\bar{\gamma}, \mathcal{D}] \subset  \mathcal{D}, \quad \forall \ V \in \mathcal{D}.	 \label{eq:invariance-2}
\end{align}
Let  
\beq \label{eq:Vi}
{\bar{V}_i   = \left(I - \frac{\partial q}{\partial \xi}\frac{\partial p}{\partial \omega}\right)e_i}
\eeq
for ${i \in \{1, \ldots, \nu\}}$. Note that ${\bar{V}_i \in \mathcal{D}}$, since 
\beq \label{eq:LVipj}
L_{\bar{V}_i}p_j
=
\frac{\partial p_j}{\partial \omega}\bar{V}_i 
\stackrel{\eqref{eq:Vi}}{=} 
\frac{\partial p_j}{\partial \omega}\left(I - \frac{\partial q}{\partial \xi}\frac{\partial p}{\partial \omega}\right)e_i 
\stackrel{\eqref{eq:condition-p-nonlinear-0}}{=} 
0 ,
\eeq
for all ${i \in \{1, \ldots, \nu\}}$ and ${j \in \{1, \ldots, r\}}$, where $p_j(\cdot)$ is the $j$-th component of the mapping $p(\cdot)$ and ${L_{V_i}p_j}(\cdot)$ is the Lie derivative of the mapping $p_j(\cdot)$ along the vector field $V_i(\cdot)$~\cite[p.8]{isidori1995nonlinear}. Then 
${\bar{V}_i \in \mathcal{D}}$ together with~\eqref{eq:invariance-1} and~\eqref{eq:invariance-2}   imply 
\begin{align}
&[\bar{\varphi}, V_i] \in  \mathcal{D}, \quad {i \in \{1, \ldots, \nu\}}, \label{eq:invariance-i-1}\\
&[\bar{\gamma}, V_i] \in  \mathcal{D},  \quad {i \in \{1, \ldots, \nu\}}.	 \label{eq:invariance-i-2}
\end{align}

We now show that~\eqref{eq:invariance-i-1} and~\eqref{eq:invariance-i-2} together imply that the constraint~\eqref{eq:PDE-nonlinear-admissible} holds. First, note that~\eqref{eq:invariance-i-2} can be expressed as 
\beq \nn
0 = \frac{\partial p}{\partial \omega} [\bar{\gamma}, V_i]
\eeq
or, equivalently, as
\beq \label{eq:invariance-i-2-a}
0 = L_{[\bar{\gamma}, V_i]} p_j 
  = L_{\bar{\gamma}} L_{V_i} p_j - L_{V_i} L_{\bar{\gamma}}  p_j 
  \stackrel{\eqref{eq:vector-fields-nonlinear},\eqref{eq:LVipj}}{=} - L_{V_i} L_{\delta}  p_j
\eeq
for all ${i \in \{1, \ldots, \nu\}}$ and ${j \in \{1, \ldots, r\}}$, where the second identity follows from~\cite[Property~(iii), p.10]{isidori1995nonlinear}. Then 
\beq \label{eq:invariance-i-2-bb}
\frac{\partial L_{\delta}  p_j}{\partial \omega}  \left(I - \frac{\partial q}{\partial \xi}\frac{\partial p}{\partial \omega}\right)e_i
\stackrel{\eqref{eq:Vi}}{=} L_{V_i} L_{\delta}  p_j
\stackrel{\eqref{eq:invariance-i-2-a}}{=} 0 
\eeq
for all ${i \in \{1, \ldots, \nu\}}$ and ${j \in \{1, \ldots, r\}}$. By integration,~\eqref{eq:invariance-i-2-bb}
yields
$$
{L_{\delta} p_j(\omega) = \left. L_{\delta} p_j(w) \right|_{w = q\circ p(\omega)}}
$$
or, equivalently,
\beq  \label{eq:invariance-i-2-b}
\frac{\partial p}{\partial \omega}(\omega) \delta(\omega) 
= \left. \frac{\partial p}{\partial w}(w) \delta (w) \right|_{w = q\circ p(\omega)}.
\eeq
Now observe that~\eqref{eq:invariance-i-1} can be equivalently expressed as 
\beq
0 
= L_{[\bar\varphi, V_i]} p_j 
= L_{\bar\varphi} L_{V_i} p_j - L_{V_i} L_{\bar\varphi} p_j \nn
\stackrel{\eqref{eq:LVipj}}{=}- L_{V_i} L_{\bar\varphi} p_j 
\stackrel{\eqref{eq:Vi}}{=} - \frac{\partial L_{\bar\varphi} p_j}{\partial \omega} \left(I - \frac{\partial q}{\partial \xi}\frac{\partial p}{\partial \omega}\right)e_i \label{eq:invariance-23},
\eeq
for all ${i \in \{1, \ldots, \nu\}}$ and ${j \in \{1, \ldots, r\}}$, where the second identity follows from~\cite[Property~(iii), p.10]{isidori1995nonlinear}. By integration,~\eqref{eq:invariance-23} yields
$$
{L_{\bar\varphi} p_j(\omega) = \left. L_{\bar\varphi} p_j(w) \right|_{w = q\circ p(\omega)}},
$$
or, equivalently,
\beq  \label{eq:invariance-XXX}
\frac{\partial p}{\partial \omega}(\omega) \bar\varphi(\omega) 
= \left. \frac{\partial p}{\partial w}(w) \bar\varphi (w) \right|_{w = q\circ p(\omega)}.
\eeq
which, in turn, using~\eqref{eq:vector-fields-nonlinear} and~\eqref{eq:invariance-i-2-b} gives~\eqref{eq:PDE-nonlinear-admissible}. We conclude that the triple $(\varphi(\xi)+\gamma(\xi)u, \kappa(\xi), p(\omega))$ satisfies the constraints of the optimization problem~\eqref{eq:optimization-problem-nonlinear}.

Finally, let ${\bar \mu = \mu(\tau(\omega))}$ and ${\bar p = p(\tau(\omega))}$. By assumption (A$_q$), after simple manipulations,~\eqref{eq:optimization-problem-proof-nonlinear} can be rewritten as
\beq \label{eq:optimization-problem-proof-nonlinear-`}
\begin{array}{ll}
    \mbox{minimize}      
    & \displaystyle\sup_{\omega \in \Omega}
    		\norm{\frac{\partial \bar\mu}{\partial \omega} (\omega) \cdot \mathbb{P}(\omega)}_{2*},
\end{array}
\eeq 
in which
\beq
\!
\mathbb{P}(\omega) = \left[I-
\frac{\partial \bar p}{\partial \omega}(\omega)^{\transpose}
\left(\frac{\partial \bar p}{\partial \omega}(\omega) \frac{\partial \bar p}{\partial \omega}(\omega)^{\transpose} \right)^{-1}
\frac{\partial \bar p}{\partial \omega}(\omega)
\right].
\eeq
A direct computation shows that ${\mathbb{P}(\omega)}$ maps covector fields in $T^{*}_\omega\Omega$ onto
$\Image(d\bar p(\omega))^{\perp}$ along $\ker(d\bar p(\omega))^{\perp}$ for every ${\omega \in \Omega}$. Furthermore, the codistributions
$\Image(d\bar p(\omega))^{\perp}$ and $\ker(d\bar p(\omega))^{\perp}$ yield an orthogonal direct sum decomposition of $T^{*}_\omega\Omega$ for every ${\omega \in \Omega}$, thus showing that ${\mathbb{P}(\omega)}$ is an orthogonal projector for every ${\omega \in \Omega}$.  A least squares argument then shows that the triple $(\varphi(\xi)+\gamma(\xi)u, \kappa(\xi), p(\omega))$ minimises the objective function of the optimization problem~\eqref{eq:optimization-problem-proof-nonlinear} and, hence, that of~\eqref{eq:optimization-problem-nonlinear}, which concludes the proof.
\end{proof}

\begin{remark} \label{rem:properties-nonlinear}
In analogy with the linear case, the available degrees of freedom in the parameterization~\eqref{eq:family-nonlinear-varphi} can be used to assign specific properties to models achieving least squares moment matching. For example,
selecting a linear signal generator~\eqref{eq:system-signal-generator} and defining ${p(\omega) = P\omega}$, ${\delta(\omega) = \Delta}$ and ${q(\xi) = Q\xi}$ yields a model~\eqref{eq:family-nonlinear} described by the mappings
\beq \label{eq:family-nonlinear-varphi-linear}
\varphi(\xi) = P(S-\Delta L)Q \xi,
\, \gamma(\xi) = P\Delta , 
\, \kappa(\xi) = h(\pi(Q\xi).
\eeq
Selecting the parameters $P$ and $\Delta$ as described in Section~\ref{ssec:properties}, it is possible to preserve the dominant eigenvalues of the \textit{linearization} of the original system (and, hence, its stability properties).
Note that, by construction, the (constant) distribution ${\ker dp}$ is $(s,l)$ invariant, with ${s(\omega) = S\omega}$ and ${l(\omega) = L\omega}$. Furthermore, the mapping $\delta(\cdot)$ is such that the distribution ${\ker dp}$ is invariant under the dynamics of the system~\eqref{eq:system-generator-surrogate-nonlinear} and such that the non-resonance condition~\eqref{eq:spectrum-constraint-nonlinear} holds.  
As a result, the parameters $p(\cdot)$, $\delta(\cdot)$ and $q(\cdot)$ are admissible for the parameterization~\eqref{eq:family-nonlinear-varphi} and the model~\eqref{eq:family-nonlinear}-\eqref{eq:family-nonlinear-varphi-linear} achieves least squares moment matching by virtue of Theorem~\ref{thm:1-nonlinear}. Finally, note that the parameterization~\eqref{eq:family-nonlinear-varphi-linear} is described by a linear differential equation with a nonlinear output map. As pointed out in~\cite{astolfi2010model}, this structure has the key advantage  that the computation of (an approximation of) the model boils down to the computation of (an approximation of) the output map. This computation can be carried out in the spirit of the results in~\cite[Chapter 4]{huang2004nonlinear}. 
\end{remark}

\subsection{A geometric interpretation} \label{ssec:geometric-interpretation}

We conclude this section with a geometric interpretation of least squares moment matching for the family of models~\eqref{eq:family-nonlinear}-\eqref{eq:family-nonlinear-varphi}.

\begin{theorem} \label{thm:least-squares-is-moment-matching-along-a-manifold}
Consider system~\eqref{eq:system-nonlinear} and the family of models~\eqref{eq:family-nonlinear}-\eqref{eq:family-nonlinear-varphi}.
Let $p(\cdot)$, $\delta(\cdot)$ and $q(\cdot)$ be {admissible} for the parameterization~\eqref{eq:family-nonlinear-varphi}. Then the model~\eqref{eq:family-nonlinear}-\eqref{eq:family-nonlinear-varphi} matches the moment of system~\eqref{eq:system-nonlinear} at $(s,l)$ along the manifold 
\beq \label{eq:manifold-nonlinear}
\mathcal{M} = \left\{\,\omega\in\Omega\,:\, \omega = q(\xi), \ \xi\in \R^r \,\right\}.
\eeq 
\end{theorem}

\begin{proof}
Let $p(\cdot)$, $\delta(\cdot)$ and $q(\cdot)$ be admissible parameters for~\eqref{eq:family-nonlinear-varphi}. 
By Theorem~\ref{thm:1-nonlinear}, 
the moment of the model~\eqref{eq:family-nonlinear}-\eqref{eq:family-nonlinear-varphi} at $(s,l)$ is well-defined.
Furthermore, by assumption (A$_{p}$), the mapping $p(\cdot)$ is a surjective submersion. Thus, for every ${\omega \in \mathcal{M}}$, there exists ${\xi \in \R^r}$ such that ${\omega = q(\xi)}$. Then
\beq \nn
\kappa (p(\omega))
\!\stackrel{\eqref{eq:manifold-nonlinear}}{=}\!  \kappa(p (q (\xi)) 
\!\stackrel{\eqref{eq:condition-p-nonlinear-0}}{=}\!  \kappa(\xi) 
\!\stackrel{\eqref{eq:family-nonlinear-varphi}}{=}\!  h(\pi(q(\xi))) 
\!\stackrel{\eqref{eq:manifold-nonlinear}}{=}\!   h(\pi(\omega)),
\eeq
which shows that the model~\eqref{eq:family-nonlinear}-\eqref{eq:family-nonlinear-varphi} matches the moment of system~\eqref{eq:system-nonlinear} at $(s,l)$ along the manifold $\mathcal{M}$.
\end{proof}

\noindent
Theorem~\ref{thm:least-squares-is-moment-matching-along-a-manifold} formalizes the idea that least squares moment matching corresponds to restricting the set where the conditions~\eqref{eq:moment-matching-PDE-model} and \eqref{eq:moment-matching-nonlinear} are required to hold. In other words, achieving least squares moment matching corresponds to achieving moment matching along the manifold $\mathcal{M}$ defined in~\eqref{eq:manifold-nonlinear}. As a direct consequence, Theorem~\ref{thm:least-squares-is-moment-matching-along-a-manifold}  provides a new perspective on least squares moment matching for linear systems, which boils down to requiring that the moment matching conditions~\eqref{eq:model-condition-2} and~\eqref{eq:Sylvester-equation-astolfi-model} hold only along the subspace ${\mathcal{S}  = \Image Q}$.

\section{Examples}  \label{sec:examples}

We now illustrate our model reduction framework with two worked-out numerical examples.
Simulations\footnote{The MATLAB code is publicly available at:~\url{https://github.com/albertopadoan/LSMR}.} have been performed in double precision using standard routines of MATLAB (version 2020b) and a 3.5 GHz Intel Core i7 processor.

\subsection{Flexible space structure}
Consider the flexible space structure benchmark model from~\cite{gawronski1991model} (see also~\cite{gawronski1996balanced}). The system is described by the equations~\eqref{eq:system-linear}, with
\bseq \label{eq:system-FSS-matrices}
\begin{align}
A &= \diag(A_1,\ldots,A_{K}), \\
B &= [\, B_1^{\transpose} \ \cdots \ B_{K}^{\transpose} \,]^{\transpose}, \\
C &= [\, C_1 \ \cdots \ C_{K} \,],
\end{align}
\eseq
in which the integer ${K>0}$ is the number of modes of the structure and 
\beq \nn
A_k = 
\bma
\begin{array}{cc}
-2\chi_k \phi_k & -\phi_k \\
\phi_k & 0
\end{array}
\ema, \,
B_k = 
\bma
\begin{array}{c}
b_k  \\
0
\end{array}
\ema, \,
C_k = 
\bma
\begin{array}{c}
c_{rk}\\
\tfrac{c_{dk}}{\phi_k}
\end{array}
\ema^{\transpose},
\eeq
where ${\chi_k \in (0,0.001)}$, ${\phi_k \in (0,100)}$,  ${b_k \in (0,1)}$ and ${C_k \in (0,10)^{1\times 2}}$ are uniformly distributed random numbers (generated in MATLAB with the function \texttt{rand} and seed $1009$). The number of modes selected for the simulations is ${K=30}$ and, hence, the order of the system is ${n=60}$.

For illustration, suppose we wish to construct a reduced order model that approximates well the original system at low frequency, say below $20$ rad/s. Further, suppose the reduced order model needs to preserve the dominant eigenvalues of the original system (and, hence, its stability properties). The parameterization~\eqref{eq:family-linear} is used to build a reduced order model of order ${r=10}$ which meets the desired specifications and achieves least squares moment matching at $\{\pm \iota\omega_i\}^{12}_{i=1}$, with ${\omega_1 = 0.01}$, ${\omega_2 = 0.1}$, ${\omega_3 = 1}$, ${\omega_4 = 5.5}$, ${\omega_5 = 10}$, ${\omega_6 = 16}$, ${\omega_7 = 20}$, ${\omega_8 = 30}$, ${\omega_9 = 50}$, ${\omega_{10} = 100}$, ${\omega_{11} = 1000}$, and ${\omega_{12} = 10000}$.
To this end, we first construct a surrogate reduced order model described by the equations~\eqref{eq:system-rom-surrogate}. The matrices $S$ and $L$ are defined as 
\beq \label{eq:system-FSS-signal-generator}
S = \diag(S_1,S_2,\ldots,S_{11},S_{12}), \quad 
L = \tfrac{1}{\sqrt{24}}[ \ \underbrace{1 \  \cdots \   1}_{24} \ ], 
\eeq
where
\beq \nn
S_i = 
\bma
\begin{array}{cc}
0 & \omega_i\\
-\omega_i & 0  
\end{array}
\ema, \hfill \quad i \in \{1, \ldots, 12\}.
\eeq 
Note that~\eqref{eq:system-FSS-signal-generator} defines an observable signal generator of order ${\nu=24}$ described by the equations~\eqref{eq:system-signal-generator} such that ${\spectrum{S} = \{\pm \iota\omega_i\}^{12}_{i=1}}$.  The Sylvester equation~\eqref{eq:Sylvester-equation-astolfi} is solved (in MATLAB with the function \texttt{sylv}) and the solution $\Pi$ is used to define the matrix ${H = C\Pi}$. A standard pole placement algorithm (implemented in MATLAB by the function \texttt{place}) is used to select the vector $\Delta$ in such a way that the matrix $(S-\Delta L)$ preserves the first ${\nu = 24}$ dominant eigenvalues of the original system. 
The matrix $P$ is defined as ${P = [\, P_1^{\transpose} \ \cdots \ P_{r}^{\transpose} \,]^{\transpose},}$ with $\{P_1,\ldots, P_r\}$ a real Jordan basis of right eigenvectors of the matrix ${S-\Delta L}$ corresponding to the first ${r = 10}$ dominant eigenvalues of the surrogate reduced order model. The matrix $Q$ is defined as ${Q = P^{\dagger}}$. By the discussion in Section~\ref{ssec:properties}, this ensures that the parameters $\Delta$, $P$, and $Q$ are admissible and, hence, that the reduced order model model achieves least squares moment matching at $\{\pm \iota\omega_i\}^{12}_{i=1}$, in agreement with Theorem~\ref{thm:moments-optimization-linear}. Furthermore, this ensures that the reduced order model preserves the first ${r = 10}$ dominant eigenvalues of the original system.

\begin{figure}[t!]
\centering
\input{FSS_rom.tex}
\centering
\caption{Top: Frequency response of system~\eqref{eq:system-FSS-matrices} (solid) and of the corresponding reduced order model (dashed). Bottom: Frequency response of the corresponding relative error.
} 
\label{fig:FSS_rom} 
\end{figure}
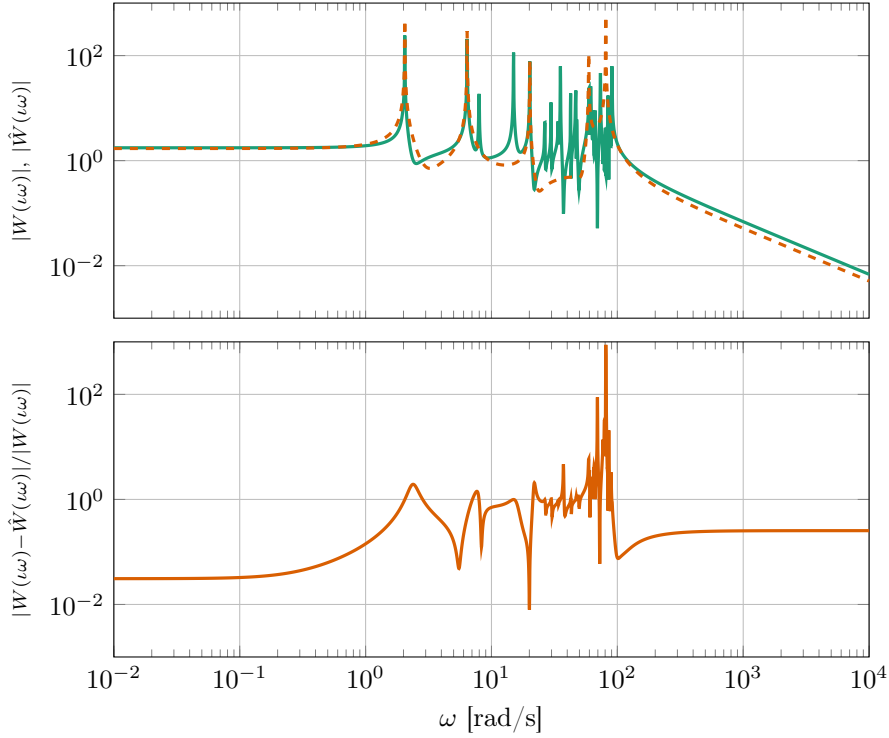%

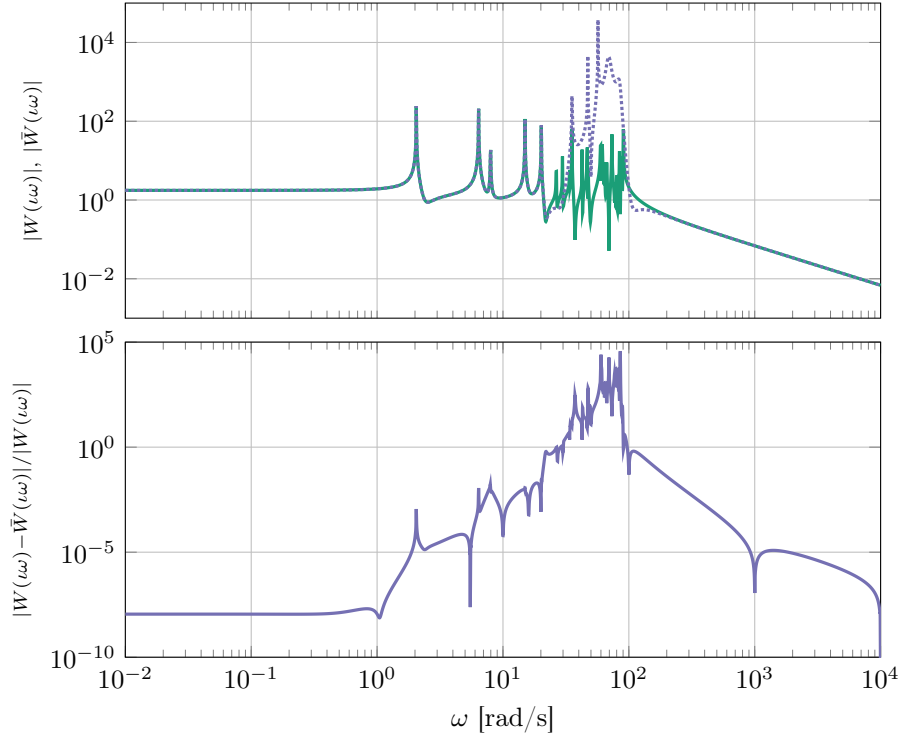
\begin{figure}[t!]
\centering
\input{FSS_rom_sur.tex}
\centering
\caption{Top: Frequency response of system~\eqref{eq:system-FSS-matrices} (solid) and of the corresponding surrogate reduced order model (dotted). Bottom: Frequency response of the corresponding relative error.
} 
\label{fig:FSS_rom_sur}
\end{figure}

The top part of Fig.~\ref{fig:FSS_rom} shows the frequency response of system~\eqref{eq:system-FSS-matrices} (solid) and of the reduced order model (dashed), respectively. The bottom part of Fig.~\ref{fig:FSS_rom} shows the frequency response of the corresponding relative error. The top part of Fig.~\ref{fig:FSS_rom_sur} shows the frequency response of system~\eqref{eq:system-FSS-matrices} (solid) and of the surrogate reduced order model (dotted), respectively. The bottom part of  Fig.~\ref{fig:FSS_rom_sur} shows the frequency response of the corresponding relative error. We observe that the reduced order model approximates relatively well the original system at low frequency (below $20$ rad/s) and preserves the first $r$ dominant eigenvalues of the original system, as required. 
Furthermore,  selecting for illustration the initial condition  ${\omega(0)=L^{\transpose}}$ yields ${\gamma_{rms} = \norm{e_{ss}}_{rms} \approx 0.1218}$  (which can be computed in MATLAB with the function \texttt{rms} using the fact that ${\norm{\omega}_{rms} = 1}$) and ${\norm{C\Pi - HP}_{2*} \approx  0.5871}$, in agreement with Theorem~\ref{thm:rms-linear}. Finally, in line with the discussion in Remark~\ref{rem:convex}, we observe that one obtains the same reduced order model by selecting the matrix $P$ as described above and by solving the optimization problem~\eqref{eq:optimization-problem-relaxed} directly using CVX~\cite{cvx2014}.

\subsection{Nonlinear inverter chain} \label{ssec:nonlinear-inverter-chain}

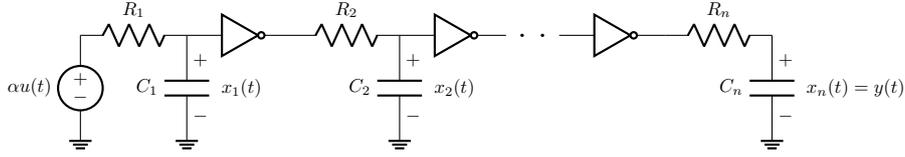
\begin{figure*}[t!]
\centering
\input{nl_inverter_chain_tikz.tex}
\centering
\caption{The chain of inverters described by system~\eqref{eq:NL_inverter_chain} in Section~\ref{ssec:nonlinear-inverter-chain}.} 
\label{fig:nl_inverter_chain}
\end{figure*}

Consider the circuit in Fig.~\ref{fig:nl_inverter_chain}. The circuit models a chain of ${n-1}$ inverters 
inspired by the benchmark model from~\cite[p.183]{gu2011model}. Following~\cite[p.242]{feynman2018feynman}, the input-output behavior of the $i$-th inverter is modelled as
\beq \nn
x_{i+1} = g_{i+1}(x_{i}) = -V_{\textrm{dd},i} \tanh\left( \frac{x_{i}}{V_T}\right),
\eeq
where ${x_{i}(t) \in \R}$ is the voltage across the $i$-th capacitor, ${V_{\textrm{dd},i} \in \Rplus}$ is the voltage supply of the $i$-th inverter, and ${V_T \in \Rplus}$ is the voltage threshold of each inverter.
The dynamics of the circuit are described by the equations
\beq \label{eq:NL_inverter_chain}
\begin{array}{rcl}
\dot{x}_1 &=& -\frac{1}{\tau_1}x_1 + \frac{\alpha}{\tau_1}u, \\
\dot{x}_2 &=& -\frac{1}{\tau_2}x_2 + \frac{1}{\tau_2}g_2(x_1), \\
		  &\vdots & 											\\
\dot{x}_n &=& -\frac{1}{\tau_n}x_n + \frac{1}{\tau_n}g_n(x_{n-1}), \\
y &=&  x_n, 
\end{array}
\eeq
where ${\alpha \in\Rplus}$ is a dimensionless parameter and  ${\tau_i \in\R}$ is the time constant of the $i$-th branch, defined as ${\tau_i = R_i C_i}$, with ${R_i \in \Rplus}$ and ${C_i \in \Rplus}$  the resistance and the capacitance of the $i$-th resistor and the $i$-th capacitor, respectively. The parameters selected for the simulations are ${n=12}$, ${V_{T} = 0.25 \,  [\,\text{V}\,]}$, ${V_{\textrm{dd},i} = \frac{1}{4(i+1)} \, [\,\text{V}\,]}$, ${\tau_i = 4(i+1) \,  [\,\text{s}\,]}$ and ${\alpha = 4}$.

For illustration, suppose we wish to construct a reduced order model that approximates well the original system, while preserving the dominant eigenvalues of the linearization of the original system (and, hence, its stability properties). The parameterization~\eqref{eq:family-nonlinear-varphi} is used to build a reduced order model of order ${r=4}$ which meets the desired specifications and achieves least squares moment matching at $(s,l)$, with 
${\tau(\omega) = \omega}$, $s(\omega) = S\omega$, $l(\omega) = L\omega$, and
\beq \label{eq:NL_inverter_chain-signal-generator}
S = \diag(S_1,S_2,\ldots,S_{4},S_{5}), \quad 
L = \tfrac{1}{\sqrt{10}}[ \ \underbrace{1 \  \cdots \   1}_{10} \ ], 
\eeq
where
\beq \nn
S_i = 
\bma
\begin{array}{cc}
0 & \omega_i\\
-\omega_i & 0  
\end{array}
\ema, \hfill \quad i \in \{1, \ldots, 5\}.
\eeq 
Note that~\eqref{eq:NL_inverter_chain-signal-generator} defines an observable linear signal generator of order ${\nu=10}$ described by the equations~\eqref{eq:system-signal-generator} such that  $\spectrum{S} = \{\pm \iota\omega_i\}^{5}_{i=1}$, with ${\omega_1 = 1}$, ${\omega_2 = 2}$, ${\omega_3 = 3}$, ${\omega_4 = 4}$ and ${\omega_5 = 5}$.

The mapping $\delta(\cdot)$ is defined as $\delta(\xi) = \Delta$, where the vector $\Delta \in \R^r$ is selected using a standard pole placement algorithm (implemented in MATLAB by the function \texttt{place}) to ensure that the matrix $(S-\Delta L)$ preserves the first ${\nu = 10}$ dominant eigenvalues of the linearization of the original system. The mapping $p(\cdot)$ is defined as  $p(\omega) = P\omega$, where the matrix ${P\in\R^{r\times\nu}}$ is defined as ${P = [\, P_1^{\transpose} \ \cdots \ P_{r}^{\transpose} \,]^{\transpose},}$ with $\{P_1,\ldots, P_r\}$ a real Jordan basis of right eigenvectors of the matrix ${S-\Delta L}$ corresponding to its first ${r = 4}$ dominant eigenvalues.
Finally, the mapping $q(\cdot)$ is defined as ${q(\xi) = P^{\dagger}\xi}$. By Remark~\ref{rem:properties-nonlinear}, the parameters $\delta(\cdot)$, $p(\cdot)$, and $q(\cdot)$ are admissible and, hence, the reduced order model model achieves least squares moment matching at $(s,l)$, in agreement with Theorem~\ref{thm:1-nonlinear}. Furthermore, in view of Remark~\ref{rem:properties-nonlinear}, this ensures that the reduced order model preserves the first ${r=4}$ dominant eigenvalues of the linearization of the original system (and, hence, its stability properties).

\begin{figure*}[t!]
\centering
\input{nl_inverter_chain.tex}
\centering
\caption{Top: Time history of the steady-state output response of system~\eqref{eq:NL_inverter_chain} and of the approximating reduced order models: $y(t)$ (solid), $\psi^{[I]}(t)$ (dotted) and $\psi^{[III]}(t)$ (dashed). Bottom: Time history of the absolute value of the corresponding steady-state output errors ${y(t) - \psi^{[I]}(t)}$ (dotted) and ${y(t) - \psi^{[III]}(t)}$ (dashed).
} 
\label{fig:nl_inverter_chain_time}
\end{figure*}
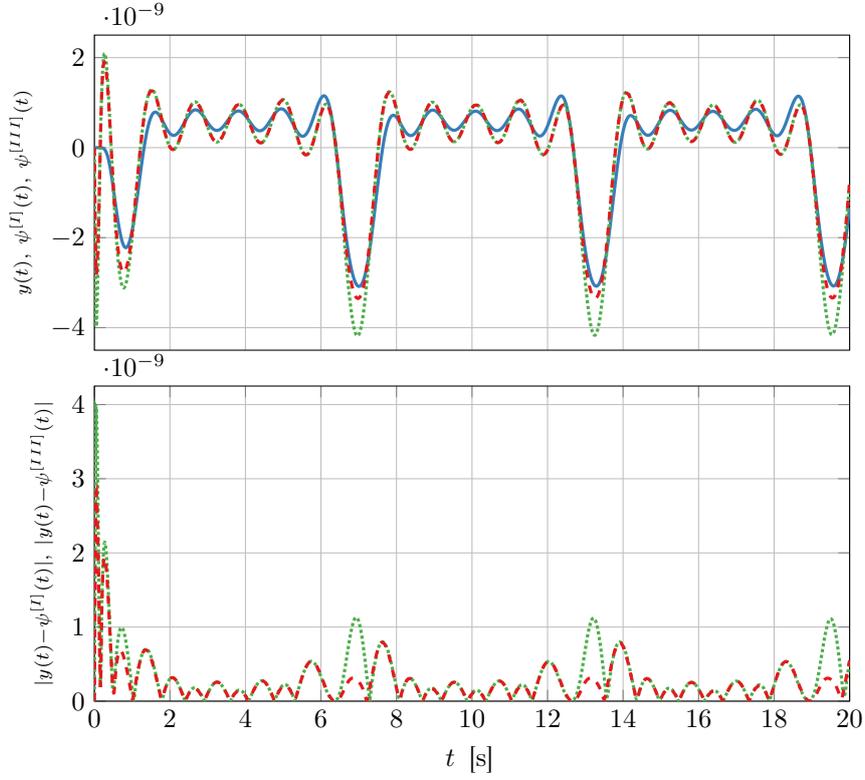

The system~\eqref{eq:NL_inverter_chain} and the corresponding reduced order model, both driven by the signal generator~\eqref{eq:system-signal-generator}-\eqref{eq:NL_inverter_chain-signal-generator}, have been numerically integrated from zero initial conditions using an explicit Runge-Kutta $(2,3)$-order integration method (implemented in MATLAB by the function \texttt{ode23}). Fig.~\ref{fig:nl_inverter_chain_time} (top) shows the output $y(t)$ of the system~\eqref{eq:NL_inverter_chain} when driven by the signal generator, and the signals $\psi^{[I]}(t)$ and $\psi^{[III]}(t)$ obtained by truncating the formal power series which defines the output $\psi(t)$ of the corresponding reduced order model to the first and third order terms, respectively. Fig.~\ref{fig:nl_inverter_chain_time} (bottom) shows the absolute value of the approximation errors ${y(t) - \psi^{[I]}(t)}$ and ${y(t) - \psi^{[III]}(t)}$. Note that, at steady-state, 
\begin{align*}
&\norm{y(t) - \psi^{[I]}(t)}_{rms} \approx 4.94\cdot 10^{-10}\\
&~>\norm{y(t) - \psi^{[III]}(t)}_{rms} \approx 3.77\cdot 10^{-10},
\end{align*}
which shows that the approximation error decreases by adding terms in the formal power series defining the output of the reduced order model.

\section{Conclusion} \label{sec:conclusion}

The model reduction problem by least squares moment matching has been studied. 
A new time-domain characterization of least squares moment matching has been introduced and used to develop a least squares model reduction theory for nonlinear systems. Families of models achieving least squares moment matching have been determined both for linear and nonlinear systems. These families have been shown to admit natural geometric and system-theoretic interpretations. The theory has been illustrated by means of numerical examples.  

A number of research directions are left open for further investigation. 
The connections between the present framework and data-driven model reduction deserve further investigation. In particular, special attention should be devoted to the relationship between least squares moment matching and the Loewner framework~\cite{mayo2007framework}. Another interesting research direction concerns the computational aspects behind least squares moment matching, including the role of weights, norms and regularization. Finally, further research is needed to clarify the connections between least squares moment matching and output regulation theory~\cite{isidori1990output}.

\section*{Acknowledgment}                               
The author warmly thanks Dr. F. Forni for suggesting the problem and Dr. A. Astolfi for his constant support.

\bibliographystyle{IEEEtran}
\bibliography{refs}

\end{document}

%% file: least-squares-moment-matching-linear-1.tex
\usetikzlibrary{patterns}
\begin{tikzpicture}[scale=0.45, every node/.style={scale=0.7}]
 
\node at (-2.5,3.5) {$\begin{array}{rcl}
\dot{ \omega} \!\!& = & \!\!  S  \omega\\
u \!\!& = & \!\!  L  \omega
\end{array}$
};
\node at (-0.5,4) {$u$};
\node at (2.3572,3.5) {
$\begin{array}{rcl}
\dot{x} & = & Ax+Bu  \\ y & = & Cx
\end{array}$};

\node at (5,4) {$y$};

\draw  (-4,4.5) rectangle (-1,2.5);
\draw  (0,4.5) rectangle (4.5,2.5);

\draw [-latex](-1,3.5) -- (0,3.5);
\draw [-latex](4.5,3.5) -- (5.5,3.5);

\node at (2.5,0.5) {
$\begin{array}{rcl}
\dot{\xi} & = & \bar{F}\xi+\bar{G}u  \\ \bar{\psi} & = & \bar{H}\xi
\end{array}$};

\node at (5,1) {$ \psi$};

\draw  (0,1.5) rectangle (4.5,-0.5);

\draw[red,thick]  (5.5,2.5) rectangle (8.5,1.5);
\node at (7,2) {$C\Pi = \bar{H}P$};
\node at (11,2) {\color{red} Moment matching};
  
\node at (7,3.5) {$C\Pi$};
\node at (11,3.5) {$A\Pi+BL=\Pi S$};
\fill[opacity=0.2,red]  (6,4) rectangle (8,3);
 
\node at (7,0.5) {$ \bar{H}P$};
\node at (11,0.5) {$ \bar{F}P+ \bar{G}L=P S$};
\fill[opacity=0.2,red]  (6,1) rectangle (8,0);

\draw [-latex](-0.5,3.5) -- (-0.5,0.5) -- (0,0.5);
\draw [-latex](4.5,0.5) -- (5.5,0.5);

\node at (-5,2) {\large  \textsc{I}};
\node at (-5,-4) {\large  \textsc{II}};

\node at (1.5,-2.5) {
$\begin{array}{rcl}
\dot{\xi} & = & F\xi+Gu  \\ \psi & = & H\xi
\end{array}$
}; 
\draw  (-1,-1.5) rectangle (4,-3.5);

\draw [-latex](-2,-2.5) -- (-1,-2.5);
\draw [-latex](4,-2.5) -- (5,-2.5);
 
\draw [ thick](-5.4284,-1) -- (13.5,-1);
\draw [ thick](-4.5,-7) -- (-4.5,5);
 
\node at (9.5,-1.6) {\color{blue} Petrov-Galerkin projection};
\node at (9.5,-2.6) {\color{blue} $P$ and $Q$ admissible};
 
\fill[opacity=0.05,blue]  (-5.4284,-1) rectangle (13.5,-7);
\fill[opacity=0.05,red]  (-5.4284,-1) rectangle (13.5,5);
 
\fill[opacity=0.15,blue]  (-2.5,-4) rectangle (-1.5,-5);
\node at (-2,-4.5) {$F$};
\node at (-1.2668,-4.5) {$=$};
\fill[opacity=0.15,blue]  (-1,-4) rectangle (0.8666,-5);
\node at (-0.0667,-4.5) {$P$};
\fill[opacity=0.4,blue]  (1,-4) rectangle (2.8666,-6);
\node at (1.9333,-5) {$\bar F$};
\fill[opacity=0.15,blue]  (3,-4) rectangle (3.8666,-6);
\node at (3.4333,-5) {$Q$};
 
\fill[opacity=0.15,blue]  (4.5,-4) rectangle (5,-5);
\node at (4.7668,-4.5) {$G$};
\node at (5.2332,-4.5) {$=$};
\fill[opacity=0.15,blue]  (5.5,-4) rectangle (7.3666,-5);
\node at (6.4333,-4.5) {$P$};
\fill[opacity=0.4,blue]  (7.5,-4) rectangle (7.9333,-6);
\node at (7.7332,-5) {$\bar G$};
 
\fill[opacity=0.15,blue]  (8.5,-4) rectangle (9.5,-4.5);
\node at (9,-4.2332) {$H$};
\node at (9.7332,-4.2999) {$=$};

\fill[opacity=0.4,blue]  (10.0667,-4.0336) rectangle (11.9333,-4.5);
\node at (11,-4.2332) {$\bar H$};
\fill[opacity=0.15,blue]  (12.0667,-4.0336) rectangle (12.9333,-6.0336);
\node at (12.5,-5.0336) {$Q$};

\draw [blue, thick] (6,-2.1) rectangle (13,-3);
\draw [dotted, thick](-1.795,-6.1);
 
\node at (-5.5,5) {};
\node at (14,5) {};
\node at (14,-7) {};
\node at (-5.5,-7) {};

\end{tikzpicture}

%% file: least-squares-moment-matching-linear-2.tex
\usetikzlibrary{patterns}
\begin{tikzpicture}[scale=0.45, every node/.style={scale=0.7}]

% Sopra

\node at (-2.5,-3) {$\begin{array}{rcl}
\dot{\bar \omega} \!\!& = & \!\! \bar S \bar \omega\\
u \!\!& = & \!\! \bar L \bar \omega
\end{array}$
};
\node at (-0.5,-2.5) {$u$};
\node at (2.3572,-3) {
$\begin{array}{rcl}
\dot{x} & = & Ax+Bu  \\ y & = & Cx
\end{array}$};

\node at (5,-2.5) {$y$};

\draw  (-4,-2) rectangle (-1,-4);
\draw  (0,-2) rectangle (4.5,-4);

\draw [-latex](-1,-3) -- (0,-3);
\draw [-latex](4.5,-3) -- (5.5,-3);

% Sotto

\node at (2.5,-6) {
$\begin{array}{rcl}
\dot{ \xi} & = &  F \xi+ Gu  \\ \psi & = &  H \xi
\end{array}$};

\node at (5,-5.5) {$ \psi$};

\draw  (0,-5) rectangle (4.5,-7);

% Additional

\draw[red,thick]  (5.5,-4) rectangle (8.5,-5);
\node at (7,-4.5) {$C\Pi = HP$};
\node at (11,-4.5) {\color{red} Moment matching};
 
%sopra
\node at (7,-3) {$C\Pi$};
\node at (11,-3) {$A\Pi+B \bar L=\Pi \bar S$};
\fill[opacity=0.2,red]  (6,-2.5) rectangle (8,-3.5);
%sotto
\node at (7,-6) {$ HP$};
\node at (11,-6) {$ FP+ G \bar L=P \bar S$};
\fill[opacity=0.2,red]  (6,-5.5) rectangle (8,-6.5);

\draw [-latex](-0.5,-3) -- (-0.5,-6) -- (0,-6);
\draw [-latex](4.5,-6) -- (5.5,-6);

% Step I
\node at (-5,2) {\large  \textsc{I}$^*$};
\node at (-5,-4) {\large  \textsc{II}$^*$};

\node at (2.5,3.5) {
$\begin{array}{rcl}
\dot{ \omega} \!\!& = & \!\! \bar S \omega\\
\bar \theta \!\!& = & \!\! \bar L \omega
\end{array}$
}; 
\draw  (1,4.5) rectangle (4,2.5);

\draw [-latex](0,3.5) -- (1,3.5);
\draw [-latex](4,3.5) -- (5,3.5);

% separation
\draw [ thick](-5.4284,-1) -- (13.5,-1);
\draw [ thick](-4.5,-7) -- (-4.5,5);

% text
\node at (9.5,4.5) {\color{blue}  Petrov-Galerkin projection};
\node at (9.5,3.5) {\color{blue} $P$ and $Q$  admissible};

%% Shades

\fill[opacity=0.05,red]  (-5.4284,-1) rectangle (13.5,-7);
\fill[opacity=0.05,blue]  (-5.4284,-1) rectangle (13.5,5);

% Model reduction by projection 

% F 
\fill[opacity=0.15,blue]  (-1.795,2) rectangle (-0.795,1);
\node at (-1.295,1.5) {$\bar S$};
\node at (-0.5618,1.5) {$=$};
\fill[opacity=0.15,blue]  (-0.295,2) rectangle (1.5716,1);
\node at (0.6383,1.5) {$P$};
\fill[opacity=0.4,blue]  (1.705,2) rectangle (3.5716,0);
\node at (2.6383,1) {$S$};
\fill[opacity=0.15,blue]  (3.705,2) rectangle (4.5716,0);
\node at (4.1383,1) {$Q$};
% % omega_0
% \fill[opacity=0.15,blue]  (4.5,8) rectangle (5,7);
% \node at (4.7668,7.5) {$\bar \omega_0$};
% \node at (5.2332,7.5) {$=$};
% \fill[opacity=0.15,blue]  (5.5,8) rectangle (7.3666,7);
% \node at (6.4333,7.5) {$P$};
% \fill[opacity=0.4,blue]  (7.5,8) rectangle (7.9333,6);
% \node at (7.7332,7) {$\omega_0$};
% L
\fill[opacity=0.15,blue]  (6.0716,2) rectangle (7.0716,1.5);
\node at (6.5716,1.7668) {$\bar L$};
\node at (7.3048,1.7001) {$=$};

\fill[opacity=0.4,blue]  (7.6383,1.9664) rectangle (9.5049,1.5);
\node at (8.5716,1.7668) {$ L$};
\fill[opacity=0.15,blue]  (9.6383,1.9664) rectangle (10.5049,-0.0336);
\node at (10.0716,0.9664) {$Q$};

\draw [blue, thick] (6,4) rectangle (13,3);
\draw [dotted, thick](-1.795,0);

%%%%%% Pointers
\node at (-5.5,5) {};
\node at (14,5) {};
\node at (14,-7) {};
\node at (-5.5,-7) {};

\end{tikzpicture}

%% file: FSS_rom.tex
\definecolor{col1}{RGB}{217,95,2}
\definecolor{col2}{RGB}{117,112,179}
\definecolor{col3}{RGB}{27,158,119} 
\pgfplotsset{compat=1.14}
\begin{tikzpicture}
\matrix{
\begin{loglogaxis}[
height = 0.475\textwidth,
width = 0.95\textwidth,
xmin=10^-2,xmax=10^4,
ymin=10^-3,ymax=10^3,
xmajorticks=false,
ylabel={{$ \scriptstyle  |W(\iota\omega)|, \ |\hat{W}(\iota\omega)|$}},
grid
]
\addplot [col3, very thick, smooth]  table [x index = {0}, y index = {1}, col sep=comma]{FSS_sys.csv};
\addplot [col1, very thick, smooth, dashed] table [x index = {0}, y index = {1}, col sep=comma]{FSS_rom.csv};
\end{loglogaxis}  
\\[0.3cm]
\begin{loglogaxis}[
height = 0.475\textwidth,
width = 0.95\textwidth,
xmin=10^-2,xmax=10^4,
ymin=10^-3,ymax=10^3,
xlabel={{$ \omega ~\text{[rad/s]} $}},
ylabel={{$ \scriptstyle  |W(\iota\omega) - \hat{W}(\iota\omega)|/|W(\iota\omega)|$}},
grid
]
\addplot [col1, very thick, smooth]  table [x index = {0}, y index = {1}, col sep=comma]{FSS_err_rom.csv};
\end{loglogaxis} 
\\[0cm]
};
\end{tikzpicture}

%% file: FSS_rom_sur.tex
\definecolor{col1}{RGB}{217,95,2}
\definecolor{col2}{RGB}{117,112,179}
\definecolor{col3}{RGB}{27,158,119}  
\pgfplotsset{compat=1.14}
\begin{tikzpicture}
\matrix{
\begin{loglogaxis}[
height = 0.475\textwidth,
width = 0.95\textwidth,
xmin=10^-2,xmax=10^4,
ymin=10^-3,ymax=10^5,
ylabel={{$ \scriptstyle  |W(\iota\omega)|, \ |\bar{W}(\iota\omega)|$}},
xmajorticks=false,
grid
]
\addplot [col3, very thick, smooth]  table [x index = {0}, y index = {1}, col sep=comma]{FSS_sys.csv};
\addplot [col2, very thick, smooth, densely dotted] table [x index = {0}, y index = {1}, col sep=comma]{FSS_rom_sur.csv};
\end{loglogaxis}  
\\[0.05cm]
\begin{loglogaxis}[
height = 0.475\textwidth,
width = 0.95\textwidth,
xmin=10^-2,xmax=10^4,
ymin=10^-10,ymax=10^5,
xlabel={{$ \omega ~\text{[rad/s]} $}},
ylabel={{$ \scriptstyle  |W(\iota\omega) - \bar{W}(\iota\omega)|/|W(\iota\omega)|$}},
grid
]
\addplot [col2, very thick, smooth]  table [x index = {0}, y index = {1}, col sep=comma]{FSS_err_rom_sur.csv};
\end{loglogaxis} 
\\[0cm]
};
\end{tikzpicture}

%% file: nl_inverter_chain_tikz.tex
 \usetikzlibrary{animations}
\begin{circuitikz}[american voltages, scale=0.7, every node/.style={transform shape}]
\ctikzset{logic ports=ieee,logic ports/scale=0.7}
	\draw (0,0) to [V,  l_=$\alpha u(t)$] ++(0,-2) node[tlground]{};
	\draw (0,0) to [R, l=${R_1}$] ++(2,0);			
	\draw (2,0) to [C, l_=${C_1}$, v^=${x_1(t)}$] ++(0,-2) node[tlground]{}; 
	\draw (2,0) to[inline not] ++(2,0);			
	\draw (4,0) to [R, l=${R_2}$] ++(2,0);			
	\draw (6,0) to [C, l_=${C_2}$, v^=${x_2(t)}$] ++(0,-2) node[tlground]{}; 
	\draw (6,0) to[inline not] ++(2,0);
	\draw[line width=1.25pt, line cap=round, shorten <=6, dash pattern= on 0pt off 6\pgflinewidth] (8,0) to (9,0); 
	\draw (9,0) to[inline not] ++(2,0);
	\draw (11,0) to [R, l=${R_n}$] ++(2,0);			
	\draw (13,0) to [C, l_=${C_n}$, v^=${x_n(t) = y(t)}$] ++(0,-2) node[tlground]{}; 				
\end{circuitikz}

%% file: nl_inverter_chain.tex
\definecolor{col1}{RGB}{55,126,184}
\definecolor{col2}{RGB}{77,175,74} 
\definecolor{col3}{RGB}{228,26,28} 
\pgfplotsset{compat=1.14}
\begin{tikzpicture}
\matrix{
\begin{axis}[
height = 0.475\textwidth,
width = 0.95\textwidth,
xmin=0,xmax=20,
ymin=-4.5*10^-9,ymax=2.5*10^-9,
xmajorticks=false,
ylabel={{$\scriptstyle y(t),$ $\scriptstyle \psi^{[I]}(t),$ $\scriptstyle \psi^{[III]}(t)$ }},
grid
]
\addplot[col1,very thick,smooth] table[x index={0}, y index={1},col sep=comma]{nl_inverter_chain_sys.csv};
\addplot[col2,very thick,smooth, densely dotted] table[x index={0}, y index={1},col sep=comma]{nl_inverter_chain_rom1.csv};
\addplot[col3,very thick,smooth,dashed] table[x index={0}, y index={1},col sep=comma]{nl_inverter_chain_rom3.csv};
\end{axis}  
\\[0.05cm]
\begin{axis}[
height = 0.475\textwidth,
width = 0.95\textwidth,
xmin=0,xmax=20,
ymin=0,ymax=4.25*10^-9,
xlabel={{$ t ~\text{\,[s]\,} $}},
ylabel={{$\scriptstyle |y(t)-\psi^{[I]}(t)|,$ $\scriptstyle |y(t)- \psi^{[III]}(t)|$ }},
grid
]
\addplot[col2,very thick,smooth, densely dotted  ] table[x index={0}, y index={1},col sep=comma]{nl_inverter_chain_err1.csv};
\addplot[col3,very thick,smooth,dashed] table[x index={0}, y index={1},col sep=comma]{nl_inverter_chain_err3.csv};
\end{axis} 
\\[0cm]
};
\end{tikzpicture}